\documentclass[12 pt]{article}
\counterwithin*{equation}{section}
\usepackage[all]{xy}
\usepackage{amsmath,amssymb,amsfonts,latexsym,float,graphics,epsfig}

\usepackage{setspace}
\usepackage[top=2cm, bottom=2cm, left=2cm, right=2cm]{geometry}
\usepackage[colorlinks]{hyperref}
\usepackage{tikz}
\usetikzlibrary{cd}
\newcommand{\contr}[1]{\iota_{#1}}

\newcommand{\dd}{\mathrm{d}}
\newtheorem{theorem}{Theorem}[section]

\newtheorem{corollary}[theorem]{Corollary}

\newtheorem{remark}{Remark}

\begin{document}
\title{Implicit Contact Dynamics and Hamilton-Jacobi Theory} \maketitle

\begin{center}
O\u{g}ul Esen\footnote{E-mail: 
\href{mailto:oesen@gtu.edu.tr}{oesen@gtu.edu.tr}}\\
Department of Mathematics, \\ Gebze Technical University, 41400 Gebze,
Kocaeli, Turkey.

\bigskip

\author{Manuel Lainz Valcázar\footnote
{E-mail: 
\href{mailto:manuel.lainz@icmat.es}{manuel.lainz@icmat.es}}
\\ 
Instituto de Ciencias Matematicas, Campus Cantoblanco \\ 
Consejo Superior de Investigaciones Cient\'ificas
 \\
C/ Nicol\'as Cabrera, 13--15, 28049, Madrid, Spain
}

\bigskip

Manuel de Le\'on\footnote{E-mail: \href{mailto:mdeleon@icmat.es}{mdeleon@icmat.es}}
\\ Instituto de Ciencias Matem\'aticas, Campus Cantoblanco \\
 Consejo Superior de Investigaciones Cient\'ificas
 \\
C/ Nicol\'as Cabrera, 13--15, 28049, Madrid, Spain
\\
and
\\
Real Academia Espa{\~n}ola de las Ciencias.
\\
C/ Valverde, 22, 28004 Madrid, Spain.
  
\bigskip

Cristina Sard\'on\footnote{E-mail: \href{mailto:mariacristina.sardon@upm.es}{mariacristina.sardon@upm.es}}
\\ Department of Applied Mathematics 
\\ Universidad Polit\'ecnica de Madrid 
\\ C/ Jos\'e Guti\'errez Abascal, 2, 28006, Madrid. Spain.

\end{center}

\date{ }

\bigskip

\begin{abstract}

\textit{In this paper we propose a Hamilton-Jacobi theory for implicit contact Hamiltonian systems in two different ways. One is the understanding of implicit contact Hamiltonian dynamics as a Legendrian submanifold of the tangent contact space, and another is as a Lagrangian submanifold of a certain symplectic space embedded into the tangent contact space. In these two scenarios we propose a Hamilton-Jacobi theory specifically derived with the aid of Herglotz Lagrangian dynamics generated by non-regular Lagrangian functions.} 
\smallskip


\smallskip


\end{abstract}

\tableofcontents
\setlength{\parskip}{4mm}

\onehalfspacing

\section{Introduction}

The goal of this paper is two-fold: one is to introduce implicit dynamical systems on contact manifolds and the other one is to develop a geometric Hamilton-Jacobi theory for these implicit contact systems. For this, in this introductory section we shall briefly summarize the fundamentals on the Hamilton-Jacobi theory for
implicit Hamiltonian dynamics in the symplectic framework (previously devised in \cite{EsLeSa18,EsLeSa20} by the same authors of the present manuscript) and then we shall carry the discussion to the contact framework to restate our goals more technically.    
 So, let us start by recalling some fundamentals of implicit dynamics on symplectic manifolds.

\subsection{Implicit Dynamics on Symplectic Manifolds and the Hamilton--Jacobi Theory}

Given a smooth manifold $\mathcal{Q}$, we refer to its cotangent bundle as $T^{*}\mathcal{Q}$ and it is known to admit a canonical (Liouville) one-form $\theta_\mathcal{Q}$ defined by
\begin{equation}\label{Liouville}
 \theta_\mathcal{Q} (X) =\langle T\pi_\mathcal{Q}(X),\tau_{T^*Q}(X) \rangle,
\end{equation}
where the brackets denote the dual pairing between $T^*\mathcal{Q}$ and $T\mathcal{Q}$, \cite{abraham1978foundations,arnold1989mathematical,
LeRo,LiMa87}. Here, $T\pi_\mathcal{Q}$ is the tangent mapping of the cotangent bundle projection $\pi_\mathcal{Q}:T^{*}\mathcal{Q}\mapsto \mathcal{Q}$ whereas $\tau_{T^*\mathcal{Q}}$ is the tangent bundle projection from $TT^*\mathcal{Q}$ to $T^*\mathcal{Q}$. The canonical symplectic two-form on $T^*\mathcal{Q}$ is the negative exterior derivative of the canonical one form, that is $\Omega_\mathcal{Q}=-d\theta_\mathcal{Q}$.

A Hamiltonian system on $T^*\mathcal{Q}$ is determined by the triple $\left( T^{\ast }\mathcal{Q},\Omega _{\mathcal{Q}},H\right) $, where $H$ is a Hamiltonian function. The Hamilton equations  determined by $H$ are the integral curves defined by the vector field provided by the following equation 
\begin{equation}\label{geomHeq}
 \iota_{X_H}\Omega_\mathcal{Q}=dH;
\end{equation}
$X_H$ is the Hamiltonian vector field, and $\iota_{X_H}$ is the interior derivative, that is the contraction mapping. In Darboux coordinates $(q^i,p_i)$ on $T^{*}\mathcal{Q}$, the canonical one-form reads $\theta_\mathcal{Q}=p_idq^i$ whereas the canonical symplectic two-form becomes $\Omega_\mathcal{Q}=dq^i\wedge dp_i$. In this local picture, the Hamilton equations turn out to be 
\begin{equation} \label{HamEqLoc}
\dot{q}^i=\frac{\partial H}{\partial p_i},\qquad \dot{p}_i=-\frac{\partial H}{\partial q^i}.
\end{equation}

\textbf{Geometric Hamilton-Jacobi Theory.}
Consider a Hamiltonian system $\left( T^{\ast }\mathcal{Q},\Omega _{\mathcal{Q}},H\right)$ and a one-form $\gamma$ on $\mathcal{Q}$, and define a vector field $X_H^{\gamma}$ on $\mathcal{Q}$ by
\begin{equation}\label{gammarelated-}
 X_H^{\gamma}:=T\pi_\mathcal{Q}\circ X_H\circ \gamma.
\end{equation}
where $X_H$ is the Hamiltonian vector field. 
With this, the geometric Hamilton-Jacobi theorem reads \cite{carinena2006geometric,de2010linearalmost}:
\begin{theorem} \label{HJT}
For a closed one-form $\gamma$ on $\mathcal{Q}$,  the following two conditions are equivalent:
\begin{enumerate}
\item The vector fields $X_{H}$ and $X_{H}^{\gamma }$ are $\gamma$-related, that is
\begin{equation}
T\gamma(X_H^{\gamma})=X_H\circ\gamma.
\end{equation}
\item The equation $d\left( H\circ \gamma \right)=0$ holds.
\end{enumerate}
\end{theorem}

The one-form $\gamma$ is assumed to be a closed one-form, so that by Poincar\'{e} lemma, we can say that there exists (locally) a function $W$ satisfying $dW=\gamma$. This function $W$ is called the characteristic function. One can easily see that substituting this into the condition $d\left( H\circ \gamma \right)=0$ results in the classical formulation of the Hamilton-Jacobi problem  
 \begin{equation}\label{hje}
H\left(q^{i},\frac{\partial W}{\partial q^{i}}\right)=c.
 \end{equation}
where $c$ is a constant of integration.
This realization of the Hamilton-Jacobi theory (henceforth also referred to as HJ theory) has been devised in various geometric frameworks, we cite \cite{EsLeSa18} and references therein for an extensive list of areas of application. 

\textbf{Implicit Hamiltonian Dynamics.}
In a series of papers \cite{marmo1992symmetries,mendella1995integrability,marmo1997constrained} the
geometry of implicit differential equations started to be investigated. In these works, for a dynamical system whose configuration space is a manifold $\mathcal{M}$, the implicit system is defined to be a submanifold of the tangent bundle $T\mathcal{M}$. Evidently, the integrability of such systems is not guaranteed. In \cite{mendella1995integrability}, there exists a geometric pathway to investigate the integrability of an implicit system. 

In the case of the canonical symplectic manifold $T^*\mathcal{Q}$, the image of a Hamiltonian vector field $X_H$ is a Lagrangian submanifold of the Tulczyjew's symplectic space $TT^*\mathcal{Q}$. Accordingly, 
an implicit Hamiltonian system is defined to be a (possibly non-horizontal) Lagrangian submanifold of $TT^*\mathcal{Q}$. Evidently, the dynamics corresponding to a non-horizontal Lagrangian submanifold is lacking (even locally) a Hamiltonian function. According to the Maslov-Hörmander theorem and the existence of a certain special symplectic structure over the cotangent bundle (see Appendix \ref{appendix-1} for the details),  one introduces a Morse family, that is a real-valued function $F$ on the total space of a fibration $\mathcal{P}\mapsto T^*\mathcal{Q}$ to accommodate the role of a Hamiltonian function.  Eventually, in the induced coordinates $(q^i,p_i,\lambda^a)$ on $\mathcal{P}$, the Lagrangian submanifold generated by a Morse family $F=F(q^i,p_i,\lambda^a)$   is computed to be
\begin{equation} \label{EandF}
\mathcal{B}=\left \{\left(q^i,p_i;\frac{\partial F}{\partial p_i},-\frac{\partial F}{\partial q^i}\right )\in TT^*\mathcal{Q}:\frac{\partial F}{\partial \lambda^a}=0\right\}.
\end{equation}
Here, $\lambda^a$'s are called the Lagrange multipliers. 
The implicit Hamiltonian dynamics is given by the system of implicit equations
\begin{equation} \label{HamEqLoc-imp}
\dot{q}^i=\frac{\partial F}{\partial p_i},\qquad \dot{p}_i=-\frac{\partial F}{\partial q^i},\qquad \frac{\partial F}{\partial \lambda^a}=0.
\end{equation}

\textbf{Geometric HJ Theory for Implicit Dynamics.} The authors of this manuscript have previously discussed a Hamilton-Jacobi theory for implicit Hamiltonian dynamics in \cite{EsLeSa18,EsLeSa20}. So, we briefly here recall the formulation that we used. Two ingredients of the classical Hamilton-Jacobi theory in Theorem \ref{HJT} are missing in the implicit picture: one is the lack of a Hamiltonian vector field $X_H$ and consequently the second is the lack of the reduced Hamiltonian vector field $X_H^\gamma$. As exhibited in \eqref{EandF}, nonetheless, we have a Lagrangian submanifold $\mathcal{B}$ instead of a Hamiltonian vector field. Accordingly, we accommodate the role of the reduced dynamics $X_H^\gamma$ by a submanifold of $T\mathcal{Q}$ as follows. Consider a differential one-form $\gamma$ on the base manifold $\mathcal{Q}$. Consider the restriction of the submanifold $\mathcal{B}$ to the image of $\gamma$ that is  
\begin{equation} \label{EandF-gamma}
\mathcal{B}\vert_{im({\gamma})}=\left \{\left(q^i,\gamma_i(q);\left.\frac{\partial F}{\partial p_i}\right\vert_{im({\gamma})},-\left.\frac{\partial F}{\partial q^i}\right\vert_{im({\gamma})}\right )\in TT^*\mathcal{Q}:\left.\frac{\partial F}{\partial \lambda^a}\right\vert_{im({\gamma})}=0\right\}.
\end{equation}
The submanifold $\mathcal{B}\vert_{im({\gamma})}$ exhibited in (\ref{EandF-gamma}) does not depend on the momentum variables. This enables us to project it to a submanifold $\mathcal{B}^{\gamma}$ of $T\mathcal{Q}$ by the tangent mapping $T\pi_Q$ as follows
\begin{equation} \label{E-gamma}
\mathcal{B}^{\gamma}=T\pi_\mathcal{Q}\circ \mathcal{B}\vert_{im ({\gamma})}=\left \{\left(q^i,\left.\frac{\partial F}{\partial p_i}\right\vert_{im({\gamma})}\right )\in T\mathcal{Q}:\left.\frac{\partial F}{\partial \lambda^a}\right\vert_{im({\gamma})}=0\right\}.
\end{equation}
Note that being a non-horizontal submanifold, $\mathcal{B}^{\gamma}$ defines an implicit differential equation on $\mathcal{Q}$. We state the generalization of the Hamilton-Jacobi Theorem \ref{HJT} as follows. 
\begin{theorem} \label{nHJT}
The following conditions are equivalent for a closed one-form $\gamma$: 
\begin{enumerate}
\item The Lagrangian submanifold $\mathcal{B}$ in (\ref{EandF}) and the submanifold $\mathcal{B}^{\gamma}$ in (\ref{E-gamma}) are $\gamma$-related, that is
$T\gamma(\mathcal{B}^{\gamma})=\mathcal{B}\vert_{im (\gamma)}$.

\item The equation $dF(q,\gamma(q),\lambda) =0$ is fulfilled, where $F$ is the Morse family generating $\mathcal{B}$.
\end{enumerate}
\end{theorem}
Assume now that the one-form $\gamma$ is exact so that $\gamma=dW$ for some real valued function $W$. Then the second condition in Theorem \ref{nHJT} gives the implicit Hamilton-Jacobi (IHJ) equations
 \begin{equation}\label{iHJEq}
   F\left(q,\frac{\partial W}{\partial q}, \lambda\right) =c.
 \end{equation} 
 where $c$ is also a constant. In Theorem  \ref{nHJT}, if the Lagrangian submanifold $\mathcal{B}$ in (\ref{EandF}) is the image of a Hamiltonian vector field $X_H$, then $\mathcal{B}^\gamma$ in (\ref{E-gamma}) becomes the image space of the vector field $X_H^\gamma$ then we retrieve the classical Hamilton-Jacobi Theorem \ref{HJT}.

\subsection{Implicit Dynamics on Contact Manifolds and the Hamilton--Jacobi Theory}\label{sub-intro}

Consider the fibration $\mathcal{Q}\mapsto \mathbb{R}$. Its first jet bundle can be identified with the extended cotangent bundle $\mathcal{T}^*\mathcal{Q}$ that is the product space $T^*\mathcal{Q}\times \mathbb{R}$. On $\mathcal{T}^*\mathcal{Q}$, there is a contact one-form 
\begin{equation}\label{eta-Q}
\eta_{\mathcal{Q}}:=dz-\theta_{\mathcal{Q}}
\end{equation}
  where  $z$ is the coordinate on the fiber $\mathbb{R}$, and $\theta_\mathcal{Q}$ is the (pull-back of the) canonical one-form in \eqref{Liouville}. Notice that we have employed abuse of notation by identifying $z$ and $\theta_{\mathcal{Q}}$ with their pull-backs on the total space $\mathcal{T}^*\mathcal{Q}$ (in Subsection \ref{Con-Man-Sec}, we will recall contact manifolds with more details). A submanifold of $\mathcal{T}^*\mathcal{Q}$ is said to be Legendrian if it is maximally isotropic. Our fist novel result in this work is Theorem \ref{thm-gen-poin-con} stating the existence of a Morse family for a given Legendrian submanifold. This result is crucial when we are analyzing  implicit dynamical equations. 
For a given Hamiltonian function $H$ on $\mathcal{T}^*\mathcal{Q}$, there are mainly two types of dynamical systems, namely contact Hamiltonian and evolution Hamiltonian formulations. 

\textbf{Contact Hamiltonian Dynamics.} 
On the contact manifold $(\mathcal{T}^*\mathcal{Q},\eta_{\mathcal{Q}}$), the contact Hamiltonian dynamics (from now on $X_H$ will denote a contact Hamiltonian vector field) is defined to be
\begin{equation}\label{contact-intro}
\iota_{X_{H}}\eta_{\mathcal{Q}} =-H,\qquad \iota_{X_{H}}d\eta_{\mathcal{Q}} =dH-\mathcal{R}(H) \eta_{\mathcal{Q}},   
\end{equation}%
where $\mathcal{R}=\partial/\partial z$ is the Reeb vector field. We cite some recent studies concerning contact Hamiltonian dynamics \cite{Br17,BrCrTa17,LeLa19,de2020review,GaGrMuMiRo20}. The dissipative nature of the contact Hamiltonian dynamics makes it proper as a geometric background for describing dissipative dynamical systems, especially irreversible processes such as thermodynamical systems \cite{bravetti2019contact,mrugala1991contact,grmela2014contact,
grmela2021multiscale,PavelkaKlikaGrmela2018,simoes2020contact}.    In Darboux coordinates $(q^i,p_i,z)$ on $\mathcal{T}^*\mathcal{Q}$, the dynamics generated by a contact Hamiltonian vector field $X_H$ is computed to be 
 \begin{equation}\label{conham-intro}
	\frac{d q^i}{dt}= \frac{\partial H}{\partial p_i}, \qquad \frac{d p_i}{dt}  = -\frac{\partial H}{\partial q^i}- 
	p_i\frac{\partial H}{\partial z}, \quad \frac{d z}{dt} = p_i\frac{\partial H}{\partial p_i} - H.
\end{equation}
This is a system of explicit differential equations. 
The extended tangent bundle $\mathcal{T}\mathcal{T}^*\mathcal{Q}$, that is the product space $T\mathcal{T}^*\mathcal{Q}\times \mathbb{R}$, admits a contact structure defined by means of a linear pencil of the vertical and the complete lifts of the canonical contact structure $\eta_{\mathcal{Q}}$ on $\mathcal{T}^*\mathcal{Q}$. Theorem \ref{the-tan-con} determines the lifted contact one-form, also known as tangent contact structure to the extended tangent bundle of the canonical contact manifold, see also \cite{IbLeMaMa97}. 
In this context, the image of a pair $(X_{H},\mathcal{R}(H))$, consisting of a contact Hamiltonian vector field $X_{H}$ and the directional derivative $\mathcal{R}(H)$ of the Hamiltonian function along the Reeb field $\mathcal{R}$, is a (horizontal) Legendrian submanifold of $\mathcal{T}\mathcal{T}^*\mathcal{Q}$. These formulations will be given more explicitly in Subsection \ref{Sec-Ham-as-Leg}. 

  \textbf{Evolution Hamiltonian Dynamics.} 
The evolution Hamiltonian dynamics $\varepsilon_H$ on a contact manifold is defined by the following equations \cite{simoes2020contact,simoes2020geometry}
\begin{equation}\label{evo-def-intro} 
 \iota_{\varepsilon_H} \eta_{\mathcal{Q}}=0, \qquad
\mathcal{L}_{\varepsilon_H}\eta_{\mathcal{Q}}=dH-\mathcal{R}(H)\eta_{\mathcal{Q}},
\end{equation}  
where $\mathcal{L}$ is the Lie derivative.
 In Darboux' coordinates $(q^i,p_i,z)$ on $\mathcal{T}^*\mathcal{Q}$, the evolution Hamiltonian dynamics is a system of differential equations
  \begin{equation}\label{conham-evo-intro}
	\frac{d q^i}{dt}= \frac{\partial H}{\partial p_i}, \qquad \frac{d p_i}{dt}  = -\frac{\partial H}{\partial q^i}- 
	p_i\frac{\partial H}{\partial z}, \quad \frac{d z}{dt} = p_i\frac{\partial H}{\partial p_i}.
\end{equation}

The kernel of the contact one-form $\eta_\mathcal{Q}$ determines a symplectic subbundle $\mathfrak{H}\mathcal{T}^*\mathcal{Q}$ of the tangent bundle $T\mathcal{T}^*\mathcal{Q}$. Moreover, the product space  $\mathfrak{H}\mathcal{T}^*\mathcal{Q} \times \mathbb{R}$ turns out to be a symplectic manifold, \cite{esen2021contact}. This is stated in Theorem \ref{the-symp} in Subsection \ref{Sec-Ev-asLag}. 
The image space of a pair $(\varepsilon_{H},\mathcal{R}(H))$, where $\mathcal{R}(H)$ is an evolution Hamiltonian vector field, is a Lagrangian submanifold of the symplectic manifold $\mathfrak{H}\mathcal{T}^*\mathcal{Q} \times \mathbb{R}$. This is given in Corollary \ref{corola}. 

\textbf{Goal I: Implicit Contact Hamiltonian Dynamics and HJ Theorem.} Consider the tangent contact manifold $\mathcal{T}\mathcal{T}^*\mathcal{Q}$. 
There are non-horizontal Legendrian submanifolds of $\mathcal{T}\mathcal{T}^*\mathcal{Q}$ that cannot (not even locally) be the image of a vector field $X_H$. We shall call the non-horizontal Legendrian submanifolds of $\mathcal{T}\mathcal{T}^*\mathcal{Q}$ as \textit{implicit contact Hamiltonian dynamics}, since the dynamics that they determine is a system of implicit differential equations. Our first goal in this paper to introduce the notion of implicit contact Hamiltonian dynamics (see  Subsection \ref{impconHam-Sec}) and write a proper Hamilton-Jacobi theory (the geometric Hamilton-Jacobi theory for explicit  Hamiltonian contact dynamics \eqref{conham-intro} has been recently examined in \cite{de2021hamilton,LeSa17,GrPa20}). In Subsection \ref{Sub-HJ-conHam}, we shall provide a particular instance (as Theorem \ref{HJT-Con}) of the contact HJ theorem. Then, in Subsection \ref{Sub-HJ-imp-con}, we shall state a geometric Hamilton-Jacobi Theorem \ref{HJ-imp-con} appropriate for implicit contact Hamiltonian dynamics.

  \textbf{Goal II: Implicit Evolution Hamiltonian Dynamics and HJ Theory.} Evidently, there exist non-horizontal Lagrangian submanifolds of $\mathfrak{H}\mathcal{T}^*\mathcal{Q} \times \mathbb{R}$  that cannot be written (not even locally) as the image space of a pair $(\varepsilon_{H},\mathcal{R}(H))$. In Subsection \ref{Imp-Evo-Sec}, we shall call non-horizontal Lagrangian submanifolds of $\mathfrak{H}\mathcal{T}^*\mathcal{Q} \times \mathbb{R}$ as \textit{implicit evolution Hamiltonian dynamics}. As in the case of the contact Hamiltonian formulation, implicit evolution Hamiltonian dynamics is provided by a system of implicit differential equations. When we construct the HJ theory in this setting we shall provide two theorems. In Subsection \ref{Sub-HJ-evol}, Theorem \ref{HJT-Evo-Con-1} will provide a Hamilton-Jacobi theorem for  evolution Hamiltonian dynamics. In Section \ref{Sub-HJ-imp-evo-11}, Theorem \ref{nHJT1-evo} provides a generalization of this theorem to implicit evolution dynamics. 

\textbf{Goal III: Physical applications.} We provide some applications of the constructed Hamilton--Jacobi theories in Section \ref{Sec-Lag}. 
The applications consist on providing Hamilton-Jacobi theories for Herglotz Lagrangian dynamics, for both  contact and evolutionary formulations. The interest of this result is that these Hamilton-Jacobi theories will be applicable even in those generated by singular (degenerate) Lagrangian functions.

\textbf{Notation.} Consider the extended configuration space $\mathcal{Q}\times \mathbb{R}$ and take it as the total space of the standard fibration from $\mathcal{Q}\times \mathbb{R}$ to $\mathbb{R}$, where the fibration is simply the projection to the second factor. The first jet manifold is diffeomorphic to $\mathcal
{T}\mathcal{Q}=T\mathcal{Q}\times \mathbb{R}$.

There exist two projections 
\begin{equation}\label{Sec-Jet-T}
\begin{split}
\tau^1_\mathcal{Q}&:\mathcal{T}\mathcal{Q}=T \mathcal{Q}\times \mathbb{R}\longrightarrow T\mathcal{Q},\qquad (U,z)\mapsto U
\\
\tau^0_\mathcal{Q}&:\mathcal{T} \mathcal{Q}=T \mathcal{Q}\times \mathbb{R}\longrightarrow \mathcal{Q}, \qquad (U,z)\mapsto \tau_\mathcal{Q}(U),
\end{split}
\end{equation}
where $\tau_\mathcal{Q}$ is the tangent bundle projection. We will refer to this space as the extended tangent bundle.   

Consider a trivial line bundle over a manifold given by $\mathcal{Q}\times \mathbb{R}\mapsto \mathcal{Q}$. The first jet bundle, denoted by $\mathcal
{T}^*\mathcal{Q}$ is diffeomorphic to the product space $T^*\mathcal{Q}\times \mathbb{R}$ that is,
\begin{equation} \label{ext-cot}
\mathcal{T}^*\mathcal{Q}=T^*\mathcal{Q}\times \mathbb{R}.
\end{equation}
We will refer to this space as the extended cotangent bundle.
There exist two projections
\begin{equation}\label{Sec-Jet}
\begin{split}
\pi^1_\mathcal{Q}&:\mathcal{T}^*\mathcal{Q}=T^*\mathcal{Q}\times \mathbb{R}\longrightarrow T^*\mathcal{Q},\qquad (\zeta,z)\mapsto \zeta
\\
\pi^0_\mathcal{Q}&:\mathcal{T}^*\mathcal{Q}=T^*\mathcal{Q}\times \mathbb{R}\longrightarrow \mathcal{Q}, \qquad (\zeta,z)\mapsto \pi_\mathcal{Q}(\zeta),
\end{split}
\end{equation}
where $\pi_\mathcal{Q}$ is the cotangent bundle projection, whereas $z$ is the standard coordinate on $\mathbb{R}$.

\section{Hamiltonian Dynamics on Contact Manifolds}

\subsection{Contact Manifolds}\label{Con-Man-Sec}

A  $(2n+1)-$dimensional manifold $\mathcal{M}$ is called a contact manifold if it is equipped with a contact one-form $\eta$ satisfying $d\eta^n
\wedge \eta \neq 0$,  \cite{arnold1989mathematical,LiMa87}. We denote a contact manifold by the pair $(\mathcal{M},\eta)$. The Reeb vector field $\mathcal{R}$ is the unique vector field  satisfying 
\begin{equation}
\iota_{\mathcal{R}}\eta =1,\qquad \iota_{\mathcal{R}}d\eta =0.
\end{equation}
At each point of the manifold $\mathcal{M}$, the kernel of the contact form $\eta$ determines the contact structure $\mathfrak{H}\mathcal{M}$. The complement of this structure, denoted by $\mathfrak{V}\mathcal{M}$,  is determined by the kernel of the exact two-form $d\eta$. These give the following decomposition of the tangent bundle 
\begin{equation}\label{decomp-TM}
T\mathcal{M}=\mathfrak{H}\mathcal{M}\oplus \mathfrak{V}\mathcal{M}, \qquad \mathfrak{H}\mathcal{M}=\ker \eta, ~ \mathfrak{V}\mathcal{M}=\ker d\eta.
\end{equation}
Here, $\mathfrak{H}\mathcal{M}$ is a vector subbundle of rank $2n$.  
The restriction of $d\eta$ to $\mathfrak{H}\mathcal{M}$ is non-degenerate so that  $(\mathfrak{H}\mathcal{M}, d\eta)$ is a symplectic vector bundle over $\mathcal{M}$. The rank of $\mathfrak{V}\mathcal{M}$ is $1$ and it is generated by  the Reeb field $\mathcal{R}$.

\textbf{Contact Diffeomorphisms.} Let $(\mathcal{M}_1,\eta_1)$ and $(\mathcal{M}_2,\eta_2)$ be two contact manifolds. A diffeomorphism $\varphi$ from $\mathcal{M}_1$ to $\mathcal{M}_2$ is said to be a contact diffeomorphism (or contactomorphism) if it preserves the contact structures, that is, $T\varphi(\ker \eta_1)=\ker \eta_2$. In terms of contact forms, a contact diffeomorphism $\varphi$ is the one satisfying 
   \begin{equation}\label{Cont-Dif}
        \varphi^*\eta_2 = \mu \eta_1.
    \end{equation}
 where $\mu $ is a non-zero conformal factor. 
To manifest the existence of this conformal factor, a mapping $\varphi$ satisfying \eqref{Cont-Dif} is also called a conformal  contact diffeomorphism. In this context, the contact mapping is denoted by the pair $(\varphi,\mu)$. If we impose that the conformal factor $\mu$ in definition \eqref{Cont-Dif} has to be equal to one, we arrive at the conservation of the contact form
  \begin{equation}
        \varphi^*\eta_2 = \eta_1.
    \end{equation}
We call such a mapping a strict contact diffeomorphism (or quantomorphism).

\textbf{Musical Mappings.} For a contact manifold $(\mathcal{M}, \eta)$, there is a musical isomorphism $\flat$ from the tangent bundle $T\mathcal{M}$ to the cotangent bundle $T^*\mathcal{M}$ defined to be 
\begin{equation}\label{flat-map}
\flat:T\mathcal{M}\longrightarrow T^*\mathcal{M},\qquad v\mapsto \iota_vd\eta+\eta(v)\eta.
\end{equation} 
This mapping takes the Reeb field $\mathcal{R}$ to the contact one-form $\eta$. We denote the inverse of this mapping by $\sharp$. Referring to this, we define a bivector field $\Lambda$ on  $\mathcal{M}$ as
\begin{equation}\label{Lambda}
\Lambda(\alpha,\beta)=-d\eta(\sharp\alpha, \sharp \beta). 
\end{equation}
The pair $(\Lambda,-\mathcal{R})$ induces a Jacobi structure \cite{Kirillov-Local-Lie,Lichnerowicz-Jacobi,Marle-Jacobi}, as a manifestation of the equalities 
\begin{equation}
[\Lambda,\Lambda]=-2\mathcal{R}\wedge \Lambda, \qquad [\mathcal{R},\Lambda]=0,
\end{equation}
where the bracket is the Schouten–Nijenhuis bracket. We cite \cite{BrGrGr17,LeMaPa97,LiMa87,Lichnerowicz-Jacobi} for more details on the Jacobi structure associated with a contact one-form. Now, referring to the bivector field $\Lambda$ we introduce the following musical mapping 
\begin{equation}\label{Sharp-Delta}
\sharp_\Lambda: T^*\mathcal{M}\longrightarrow T\mathcal{M}, \qquad  \alpha\mapsto \Lambda(\alpha,\bullet)= \sharp \alpha - \alpha(\mathcal{R})  \mathcal{R}. 
\end{equation}
Evidently, the mapping $\sharp_\Lambda$ fails to be an isomorphism. Notice that the kernel is spanned by the contact one-form $\eta$, so that the image space of $\sharp_\Lambda$ is precisely the horizontal bundle $\mathfrak{H}\mathcal{M}$ exhibited in \eqref{Horizontal-space}.

\textbf{Submanifolds.} 
Let $(\mathcal{M},\eta)$ be a contact manifold. Recall the associated bivector field $\Lambda$ defined in \eqref{Lambda}. Consider a linear subbundle $ \Xi$ of the tangent bundle $T\mathcal{M}$ (that is, a distribution on $\mathcal{M}$). We define the contact complement of $\Xi$ as
\begin{equation}
\Xi^\perp : = \sharp_ \Lambda(\Xi^o),
\end{equation}
 where the sharp map on the right hand side is the one in   \eqref{Sharp-Delta} and   $\Xi^o$ is the annihilator of  $\Xi$.
Let $\mathcal{N}$ be a submanifold of $\mathcal{M}$. We say that $\mathcal{N}$ is \cite{LeLa19,deLeon2020infinitesimal,LiMa87}:
    \begin{itemize}
        \item \emph{Isotropic} if $T\mathcal{N}\subseteq {T\mathcal{N}}^{\perp }$,
        \item \emph{Coisotropic} if $T\mathcal{N}\supseteq {T\mathcal{N}}^{\perp }$
      ,
        \item \emph{Legendrian} if $T\mathcal{N}= {T\mathcal{N}}^{\perp }$.
    \end{itemize}


\subsection{Extended Cotangent Bundle and Legendrian Submanifolds}\label{Gen-Fam-Leg}

An  example of a contact manifold is the extended cotangent bundle $\mathcal{T}^*\mathcal{Q}=T^*\mathcal{Q}\times \mathbb{R}$ given in \eqref{ext-cot}.  
Here, the contact one-form is $\eta_{\mathcal{Q}}$ given in \eqref{eta-Q}. 
There exist Darboux coordinates $(q^i,p_i,z)$ on $\mathcal{T}^*\mathcal{Q}$, where $i$ is running from $1$ to $n$ (which is the dimension of the configuration manifold). In these coordinates, the contact one-form and the Reeb vector field are computed to be
\begin{equation}
\eta _{\mathcal{Q}}= d z -  p_i d q^i, \qquad \mathcal{R}=\frac{\partial}{\partial z},
\end{equation}
respectively. Notice that, in this realization, the horizontal bundle is generated by the vector fields
\begin{equation}\label{Horizontal-space}
\mathfrak{H}\mathcal{T}^*\mathcal{Q}=span\{\xi_i,\xi^i\},\qquad \xi_i=\frac{\partial}{\partial q^i} + p_i \frac{\partial}{\partial z},~\xi^i=\frac{\partial}{\partial p_i}.
\end{equation}
It is important to note that these generators are not closed under the Jacobi-Lie bracket, i.e., 
\begin{equation}
[\xi^{i},\xi_{j}]=\delta^i_j\mathcal{R},
\end{equation}
where $\delta^i_j$ stands for the Kronecker delta. 
 In terms of the Darboux coordinates, we compute the musical mapping $\sharp_\Lambda$
  in \eqref{Sharp-Delta}   as
\begin{equation}\label{zap}
\sharp_\Lambda:\alpha_i dq^i + \alpha^i dp^i + 
u dz\mapsto 
\alpha^i \frac{\partial}{\partial q^i}-(\alpha_i + p_i u)\frac{\partial}{\partial p_i}  + \alpha^ip_i  \frac{\partial}{\partial z}. 
\end{equation}

\textbf{Legendrian Submanifolds.}
There is a nice way to characterize Legendrian submanifolds of contact manifolds. In Darboux coordinates the picture is as follows. 
Let $(\mathcal{T}^*\mathcal{Q},\eta_{\mathcal{Q}})$ be the canonical contact manifold equipped with Darboux coordinates $(q^i,p_i,z)$ for $i=1,\dots n$. 
Consider a partition $A\cup B$ of the set of indices $(1,\dots, n)$ into two disjoint subsets, so that the Darboux coordinates on the extended cotangent bundle $\mathcal{T}^*\mathcal{Q}$ turn out to be
$(q^a,q^b,p_a,p_b,z)$ where $a\in A$ and $b\in B$. For a function $\Phi( {q}^a, {p}_b)$ of $n$ variables $ {q}^a$, $a\in A$ and $ {p}_b$, $b\in B$, the
$2n + 1$ equations
\begin{equation}\label{Leg-Sub-Loc}
 {q}^b=-\frac{\partial \Phi }{\partial {p}_b},\qquad 
 {p}_a=\frac{\partial \Phi }{\partial  {q}^a}, \qquad z=\Phi- {p}_b \frac{\partial \Phi }{\partial  {p}_b}
\end{equation}
define a Legendrian submanifold of $\mathcal{T}^*\mathcal{Q}$, see \cite{arnold1989mathematical}.  Conversely, in a neighborhood of every point, every Legendrian submanifold of
$\mathcal{T}^*\mathcal{Q}$ can be given as \eqref{Leg-Sub-Loc} for at least one
of the $2^n$ possible choices of the subset $A$. We  record this generic picture of Legendrian submanifolds as
\begin{equation}\label{Leg-Gen-1}
\mathcal{N}=\left\{(q^a,q^b,p_a,p_b,z)\in \mathcal{T}^*\mathcal{Q}: {q}^b=-\frac{\partial \Phi }{\partial {p}_b},
{p}_a=\frac{\partial \Phi }{\partial  {q}^a}, z=\Phi- {p}_b \frac{\partial \Phi }{\partial  {p}_b}  \right\}.
\end{equation}
Even though this definition is determining all the Legendrian submanifolds of $\mathcal{T}^*\mathcal{Q}$, it cannot be used to examine the dynamics, since it is not given in terms of 'sections'. Let us examine now the horizontal and non-horizontal Legendrian submanifold appropriate for our interests. 
 
 \textbf{Horizontal Legendrian Submanifolds.} 
Let $W$ be a real valued function on the base manifold $\mathcal{Q}$. Its first prolongation is a section of the bundle $\pi^0_\mathcal{Q}$ displayed in  \eqref{Sec-Jet} that is,
\begin{equation}\label{j1F}
\mathcal{T}^* W:\mathcal{Q}\longrightarrow \mathcal{T}^*\mathcal{Q}=T^*\mathcal{Q}\times \mathbb{R},\qquad (q^i)\mapsto \left(q^i,\frac{\partial W}{\partial q^i},W(q)\right).
\end{equation}
The image space of the first prolongation $\mathcal{T}^* W$ is  a Legendrian  submanifold of $\mathcal{T}^*\mathcal{Q}$. The converse of this assertion is also true, that is, if the image space of a section   of $\pi^0_\mathcal{Q}$ is a Legendrian submanifold of the extended cotangent bundle $\mathcal{T}^*\mathcal{Q}$, then it is  the first prolongation of a function $W$.  Evidently, this is a particular instance of the local realization in \eqref{Leg-Gen-1}, where the index set $B$ is empty, i.e., if the generating function depends only on the coordinates of the base manifold  $\Phi=\Phi(q^i)$. 

\textbf{Non-horizontal Legendrian Submanifolds.} A Legendrian submanifold of $\mathcal{T}^*\mathcal{Q}$ is called non-horizontal if it cannot be written as the image of a section of the fibration $ \mathcal{T}^*\mathcal{Q}\mapsto  \mathcal{Q}$. In this present work, to characterize non-horizontal Legendrian submanifolds, we shall employ Morse families (see Appendix \ref{appendix-1} for the definition of Morse families). 
Consider a Morse family $E$ defined on the total space of a smooth bundle $(\mathcal{P},\tau,\mathcal{Q})$ over the base manifold $\mathcal{Q}$. A Morse family generates a Lagrangian submanifold $\mathcal{S}$ of the cotangent bundle $T^*\mathcal{Q}$. We lift this Lagrangian submanifold to a Legendrian submanifold of the contact manifold $\mathcal{T}^*\mathcal{Q}$. To see this, consider a local system of coordinates $(q^i)$ on the base manifold $\mathcal{Q}$, and the induced coordinates  $(q^i,\epsilon^a)$ on the total space $\mathcal{P}$. Then, referring to the Darboux coordinates on $\mathcal{T}^*\mathcal{Q}$, the Legendrian submanifold $\mathcal{N}$  generated by a Morse family $E=E(q,\epsilon)$ is computed to be 
\begin{equation} \label{MFGen-C}
\mathcal{N}  =\left \{\Big(q^i,\frac{\partial E}{\partial q^i}(q,\epsilon),E(q,\epsilon)\Big )\in \mathcal{T}^*\mathcal{Q}: \frac{\partial E}{\partial \epsilon^a}=0\right \}\subset \mathcal{T}^*\mathcal{Q}.
\end{equation} 
We plot this definition in the following picture
\begin{equation} \label{Morse-Gen---1}
\xymatrix{ \mathcal{T}^*\mathcal{Q}\ar@(ul,ur)^{ \mathcal{N}} \ar[dd]_{\pi_{\mathcal{Q}}^0}&&
\mathcal{P}  \ar[dd]^{\tau}\ar[r]^{E} &\mathbb{R} \\
\\
	\mathcal{Q} \ar@{=}[rr]&& \mathcal{Q} 
}
\end{equation}  
The following theorem proves the existence of Morse families generating Legendrian submanifolds. In this sense, it provides a local realization for Legendrian submanifolds, alternative to the one in \eqref{Leg-Gen-1}. The advantages of this local realization will be realized later on when we will focus on the implicit dynamics. 
\begin{theorem}\label{thm-gen-poin-con}
For any (not necessarily horizontal) Legendrian submanifold of a contact manifold there exists a Morse family generating it.
  \end{theorem}
\textbf{Proof.}
It is interesting to connect two local realizations of Legendrian submanifolds given in \eqref{Leg-Gen-1} and \eqref{MFGen-C}. To prove the assertion, consider first the generic picture in \eqref{Leg-Gen-1}, and then generate this Legendrian submanifold by a Morse family. For this, assume the partition of the index set $(q^a,q^b)$ on the base manifold $\mathcal{Q}$ and a fiber bundle $\mathcal{P}$ with coordinates $(q^a,q^b,p_b)$. Here, the projection is simply
\begin{equation}
\tau:\mathcal{P}\longrightarrow \mathcal{Q},\qquad (q^a,q^b,p_b)\mapsto (q^a,q^b).
\end{equation} 
  Then, consider the following function on the total space
  \begin{equation}\label{Morse-Fam}
  E(q^a,q^b,p_b)=\Phi(q^a,p_b)+q^bp_b
  \end{equation}
  where $\Phi=\Phi(q^a,p_b)$ is the function determining the Legendrian submanifold \eqref{Leg-Gen-1}. It is immediate to see that the function  $E$ in \eqref{Morse-Fam} satisfies  the requirement of being a Morse family given in Appendix \ref{appendix-1}. According to the definition in \eqref{MFGen-C}, the Legendrian submanifold generated by the Morse family \eqref{Morse-Fam} is computed to be
  \begin{equation} \label{MFGen-C-}
\mathcal{N}  =\left \{\Big(q^a,q^b,\frac{\partial \Phi}{\partial q^a},p_b,\Phi(q^a,p_b)+q^bp_b \Big )\in \mathcal{T}^*\mathcal{Q}: \frac{\partial \Phi}{\partial p^b}+q^b=0 \right \}.
\end{equation} 
  A direct observation shows that the submanifold in \eqref{MFGen-C-} coincides with the submanifold in \eqref{Leg-Gen-1}.

\textbf{Special Contact Structures.} Morse families are generating Legendrian submanifolds of the extended cotangent manifolds $\mathcal{T}^*\mathcal{Q}$. To transfer such a Legendrian submanifold to a contact manifold $(\mathcal{M},\eta)$ one needs to employ a special contact structure, see \cite{esen2021contact}. Let us first introduce the notion of the special contact structure and then merge it with a Morse family. 
Let $(\mathcal{M},\eta)$ be a contact manifold and the total space of a fiber bundle $(\mathcal{M},\rho,\mathcal{Q})$. We introduce a special contact structure as a quintuple
		\begin{equation} \label{scs-5}  
		(\mathcal{M},\rho,\mathcal{Q},\eta,\Phi),
				\end{equation}
where $
		\Phi$ is a fiber preserving contact
		diffeomorphism (that is $\pi_{\mathcal{Q}}^0\circ \Phi=\rho$) from $\mathcal{M}$ to the canonical contact manifold $(\mathcal{T}^*\mathcal{Q},\eta_{\mathcal{Q}})$. Here, we have a diagram exhibiting a   special contact structure in a pictorial way
		\begin{equation} \label{scs}
		\xymatrix{ \mathcal{T}^*\mathcal{Q}\ar[ddr]_{\pi^0_\mathcal{Q}} &&\mathcal{M}
			\ar[ddl]^{\rho} \ar[ll]_{\Phi}
			\\ \\&\mathcal{Q} }  
		\end{equation}
where $\pi^0_\mathcal{Q}$ is the fibration given in \eqref{Sec-Jet}. 
The pair $(\mathcal{M},\eta)$ is said to be the underlying contact manifold of the special contact structure.

\textbf{Generating Family for Legendrian Submanifolds.}
We now merge a Morse family $E$ defined on  $(\mathcal{P},\tau,\mathcal{Q})$ and a special contact space $  		(\mathcal{M},\rho,\mathcal{Q},\eta,\Phi)$ in order to arrive at a Legendrian submanifold of $(\mathcal{M},\eta)$. For this, consider the following commutative diagram 
\begin{equation} \label{Morse-Gen--}
\xymatrix{
	\mathbb{R}& \mathcal{P} \ar[dd]^{\tau}\ar[l]^{E}& \mathcal{T}^*\mathcal{Q} \ar@(ul,ur)^{ \mathcal{N}} \ar[ddr]_{\pi_{\mathcal{Q}}^0}& &\mathcal{M}\ar@(ul,ur)^{ \mathcal{N}_{E}} \ar[ll]_{\Phi} \ar[ddl]^{\rho}\\ \\ &
	\mathcal{Q} \ar@{=}[rr]&& \mathcal{Q} 
}
\end{equation}   
Referring to the definition in \eqref{MFGen-C}, we obtain a Legendrian submanifold $\mathcal{N}$ of the jet bundle $\mathcal{T}^*\mathcal{Q}$. Then, by employing the inverse of the contact diffeomorphism $\Phi$, we arrive at a Legendrian submanifold $\mathcal{N}_{E}$ of $\mathcal{M}$.

\subsection{Dynamics on Contact Manifolds}\label{Sub-dyn-con}

For a real valued function $H$ on a contact manifold $(\mathcal{T}^*\mathcal{Q},\eta_{\mathcal{Q}})$, there corresponds a contact vector field $X_H$ defined as follows
\begin{equation}
\iota_{X_{H}}\eta_{\mathcal{Q}} =-H,\qquad \iota_{X_{H}}d\eta_{\mathcal{Q}} =dH-\mathcal{R}(H) \eta_{\mathcal{Q}},   \label{contact}
\end{equation}%
where $\mathcal{R}$ is the Reeb vector field. Here, $H$ is called the (contact) Hamiltonian function and $X_H$ is called the (contact)  Hamiltonian vector field. In terms of the musical mapping $\sharp_\Lambda$   in \eqref{Sharp-Delta} a contact Hamiltonian vector field is defined as
\begin{equation}
X_{H}=\sharp_\Lambda(dH)-H\mathcal{R},
\end{equation}
where $\mathcal{R}$ is the Reeb field. 
A direct computation determines a conformal factor for a Hamiltonian vector fields as
\begin{equation}\label{L-X-eta}
\mathcal{L}_{X _{H}}\eta_{\mathcal{Q}} =
d\iota_{X _{H}}\eta_{\mathcal{Q}}+\iota_{X _{H}}d\eta_{\mathcal{Q}}= -\mathcal{R}(H)\eta_{\mathcal{Q}}.
\end{equation}
According to \eqref{L-X-eta}, the flow of a contact Hamiltonian system preserves the contact structure, but it does not preserve neither the contact one-form nor the Hamiltonian function. Instead we obtain
\begin{equation}
{\mathcal{L}}_{X _H} \, H = - \mathcal{R}(H) H.
\end{equation}
Being a non-vanishing top-form we can consider $d\eta^n
\wedge \eta$ as a volume form on $\mathcal{T}^*\mathcal{Q}$.  
Hamiltonian motion does not preserve the volume form since
\begin{equation}
{\mathcal{L}}_{X _H}  \, (d\eta_{\mathcal{Q}}^n
\wedge \eta_{\mathcal{Q}}) = - (n+1)  \mathcal{R}(H) d\eta_{\mathcal{Q}}^n
\wedge \eta_{\mathcal{Q}}.
\end{equation}%
However, it is immediate to see that, for a nowhere vanishing Hamiltonian function $H$, the quantity $ {H}^{-(n+1)}   (d\eta_{\mathcal{Q}})^n \wedge\eta_{\mathcal{Q}}$ 
is preserved along the motion (see \cite{BrLeMaPa20}).

Referring to the Darboux coordinates $(q^i,p_i,z)$, for a Hamiltonian function $H$, the Hamiltonian vector field, determined in \eqref{contact}, is computed to be
\begin{equation}\label{con-dyn}
X _H=\frac{\partial H}{\partial p_i}\frac{\partial}{\partial q^i}  - \left(\frac{\partial H}{\partial q^i} + \frac{\partial H}{\partial z} p_i \right)
\frac{\partial}{\partial p_i} + (p_i\frac{\partial H}{\partial p_i} - H)\frac{\partial}{\partial z}.
\end{equation}
We obtain that the Hamilton equations for $H$ as
\begin{equation}\label{conham}
	\frac{d q^i}{dt}= \frac{\partial H}{\partial p_i}, \qquad \frac{d p_i}{dt}  = -\frac{\partial H}{\partial q^i}- 
	p_i\frac{\partial H}{\partial z}, \quad \frac{d z}{dt} = p_i\frac{\partial H}{\partial p_i} - H.
\end{equation}

\textbf{Evolution Hamiltonian Dynamics.}
Referring to a Hamiltonian function $H$ on the contact manifold $(\mathcal{T}^*\mathcal{Q},\eta_{\mathcal{Q}})$: the \emph{evolution vector field} of $H$~\cite{simoes2020contact}, denoted as $\varepsilon_H$, which is the one that satisfies
\begin{equation}\label{evo-def} 
\mathcal{L}_{\varepsilon_H}\eta_{\mathcal{Q}}=dH-\mathcal{R}(H)\eta_{\mathcal{Q}},\qquad \eta_{\mathcal{Q}}(\varepsilon_H)=0.
\end{equation} 
 In terms of the musical mapping $\sharp_\Lambda$   in \eqref{Sharp-Delta} a evolution Hamiltonian vector field is defined as
\begin{equation}
\varepsilon_{H}=\sharp_\Lambda(dH),
\end{equation}  Let us compare the contact Hamiltonian and evolution Hamiltonian dynamics citing \cite{simoes2020geometry}.
The evolution vector field $\varepsilon_H$ preserves the Hamiltonian function $H$, but the
contact Hamiltonian dynamics $X_{H}$ dissipates it. On the other hand, $X_{H}$ preserves the contact structure that is the kernel of $\eta_\mathcal{Q}$, but $\varepsilon_H$ does not. 
In local coordinates the evolution Hamiltonian vector field is computed to be
\begin{equation}\label{evo-dyn}
	\varepsilon_H=\frac{\partial H}{\partial p_i}\frac{\partial}{\partial q^i}  - \left (\frac{\partial H}{\partial q^i} + \frac{\partial H}{\partial z} p_i \right)
	\frac{\partial}{\partial p_i} + p_i\frac{\partial H}{\partial p_i} \frac{\partial}{\partial z},
\end{equation}
so that the integral curves satisfy the evolution equations
\begin{equation}\label{evo-eq}
	\frac{d q^i}{dt}= \frac{\partial H}{\partial p_i}, \qquad \frac{d p_i}{dt}  = -\frac{\partial H}{\partial q^i}- 
	p_i\frac{\partial H}{\partial z}, \quad \frac{d z}{dt} = p_i\frac{\partial H}{\partial p_i}.
\end{equation}

\section{Implicit Contact Hamiltonian Dynamics and HJ Theory}
\label{Sec-imp-con-Ham-HJ}

\subsection{HJ Theory for Contact Hamiltonian Dynamics}\label{Sub-HJ-conHam}

Consider the extended manifold $\mathcal{Q}\times \mathbb{R}$ and the extended cotangent bundle $\mathcal{T}^{*}\mathcal{Q}\simeq T^{*}\mathcal{Q} \times \mathbb{R}$. Define a fibration 
\begin{equation}\label{rho}
\rho:\mathcal{T}^{*}\mathcal{Q}\simeq T^{*}\mathcal{Q} \times \mathbb{R}\longrightarrow \mathcal{Q}\times \mathbb{R},\qquad (w,z)\mapsto \big(\pi_{\mathcal{Q}}(w),z\big).
\end{equation}
We are particularly interested in the sections of the fibration $\rho$
of the form 
 \begin{equation}\label{gamma}
 \gamma: \mathcal{Q}\times \mathbb{R}\longrightarrow \mathcal{T}^{*}\mathcal{Q}\simeq T^{*}\mathcal{Q} \times \mathbb{R},\qquad (q,z)\mapsto (\hat{\gamma}(q),z) 
 \end{equation}
 where $\hat{\gamma}$ is a differential form on the base manifold  $\mathcal{Q}$. 
 
Consider a section $\gamma$ in form \eqref{gamma} and a contact Hamiltonian vector field $X_H$. We define a vector field $X_{H}^{\gamma}$  on the extended space $\mathcal{Q}\times \mathbb{R}$ according to the commutation of the following diagram
\begin{equation}
\xymatrix{ \mathcal{T}^{*}\mathcal{Q}\simeq T^{*}\mathcal{Q} \times \mathbb{R}
\ar[dd]^{\rho} \ar[rrr]^{X_H}&   & &T\mathcal{T}^{*}\mathcal{Q}\ar[dd]^{T{\rho}}\\
  &  & &\\
\mathcal{Q}\times \mathbb{R} \ar@/^2pc/[uu]^{\gamma}\ar[rrr]^{X^{\gamma}_H}&  & & T(\mathcal{Q}\times \mathbb{R})}
\end{equation}
that is 
\begin{equation}\label{X-gamma-con}
X_{H}^{\gamma}:=T{\rho}\circ X_H \circ \gamma,
\end{equation}
where $T{\rho}$ is the tangent mapping of the fibration $\rho$ in \eqref{rho}. Next, we state a Hamilton-Jacobi theorem for contact Hamiltonian dynamics.

\begin{theorem} \label{HJT-Con}
For a section $\gamma$ in the form \eqref{gamma} such that its first component $\hat{\gamma}$ is closed, the following two conditions are equivalent:
\begin{enumerate}
\item The vector fields $X_{H}$ and $X^{\gamma}_H$ are $\gamma$-related, that is
\begin{equation}\label{HJT-Con-1-eq-1--}
T\gamma \circ  X_H^{\gamma} =X_H\circ \gamma
\end{equation}
where $T\gamma:T(\mathcal{Q}\times \mathbb{R})\mapsto T\mathcal{T}^{*}\mathcal{Q}$ is the tangent mapping of  the section $\gamma$. 
\item The equation 
\begin{equation}\label{HJT-Con-1-eq-2--}
d( H\circ \gamma)=  \gamma^*(\mathcal{R}(H)\eta_{\mathcal{Q}}) 
\end{equation}
is fulfilled. Here, $\mathcal{R}$ is the Reeb field. 
\end{enumerate}

\end{theorem}

\textbf{Proof.} We prove this theorem in local coordinates, so that we consider the closed form $\hat{\gamma}$ as an exact one-form that is $\hat{\gamma}=dW$ for a real valued (possibly local) function $W$ on $\mathcal{Q}$. In this case, the vector field defined in \eqref{X-gamma-con} is written as
 \begin{equation}
X^{\gamma}_H(q,z)= \left. \frac{\partial H}{\partial p_i}\right\vert_{im (\gamma)} \frac{\partial}{\partial q^i} +  \Big(\frac{\partial W}{\partial q^i}\left. \frac{\partial H}{\partial p_i}\right\vert_{im (\gamma)} - \left. H\right\vert_{im (\gamma)} \Big) \frac{\partial}{\partial z}.
   \end{equation} 
In the induced coordinates $(q^i,z,\dot{q}^i,\dot{z})$ on the tangent bundle $T(\mathcal{Q}\times \mathbb{R})$, the tangent mapping of the section $\gamma$ is computed to be  
 \begin{equation}\label{T-gamma}
 T\gamma(q^i,z,\dot{q}^i,\dot{z})= \left(q^i,\frac{\partial W}{\partial q^i},z,\dot{q}^i,\frac{\partial W}{\partial q^i \partial q^j}\dot{q}^j,\dot{z}\right).
 \end{equation}
Accordingly, the left hand side of the equation \eqref{HJT-Con-1-eq-1--} is computed to be
 \begin{equation}\label{Proof-1-}
 \begin{split}
  T\gamma \circ   X^{\gamma}_H(q,z)&=\left(q^i,\frac{\partial W}{\partial q^i},z; \left.\frac{\partial H}{\partial p_i}\right\vert_{im (\gamma)}, \frac{\partial W}{\partial q^i \partial q^j}\left. \frac{\partial H}{\partial p_j}\right\vert_{im (\gamma)}
  ,\frac{\partial W}{\partial q^i}\left. \frac{\partial H}{\partial p_i}\right\vert_{im (\gamma)} - \left. H\right\vert_{im (\gamma)}\right).
 \end{split}
  \end{equation}
In this local realization, the right hand side of the first condition \eqref{HJT-Con-1-eq-1--} turns out to be
  \begin{equation}\label{Proof-2-}
   \begin{split}
 X_H\circ \gamma(q,z)&=\Big(q^i,\frac{\partial W}{\partial q^i},z ;\left. 
  \frac{\partial H}{\partial p_i}\right\vert_{im (\gamma)},-\left. \frac{\partial H}{\partial q^i}\right\vert_{im (\gamma)}
  -\frac{\partial W}{\partial q^i}\left. \frac{\partial H}{\partial z}
 \right\vert_{im (\gamma)}, \frac{\partial W}{\partial q^i}\left. \frac{\partial H}{\partial p_i}\right\vert_{im (\gamma)} - \left. H\right\vert_{im (\gamma)}
  \Big). 
   \end{split}
    \end{equation}
    In order to satisfy the identity in \eqref{HJT-Con-1-eq-1--}, we consider the expressions \eqref{Proof-1-} and \eqref{Proof-2-}. The first three entries of these local realizations are the same. Two of the fiber variables (that coincide with the fourth and sixth entries) are also the same. See that the fifth entries of \eqref{Proof-1-} and \eqref{Proof-2-} are not equal. By considering them equal, we arrive at  the identity 
 \begin{equation}\label{evo-con-HJ-expli-}
\frac{\partial W}{\partial q^i \partial q^j}\left. \frac{\partial H}{\partial p_j}\right\vert_{im (\gamma)}+\left. \frac{\partial H}{\partial q^i}\right\vert_{im (\gamma)}
  +\frac{\partial W}{\partial q^i}\left. \frac{\partial H}{\partial z}
  \right\vert_{im (\gamma)}=0.
     \end{equation}
We will prove now that the identity \eqref{evo-con-HJ-expli-} is precisely the local form of the condition in \eqref{HJT-Con-1-eq-2--} and this would complete the proof. To this end, consider the second condition \eqref{HJT-Con-1-eq-2--} and, in the local coordinates, compute the left hand side of this identity as 
      \begin{equation} \label{Condu-1} 
         \begin{split}   
     d( H\circ \gamma) = 
     & = \left(\left. \frac{\partial H}{\partial q^i}\right\vert_{im (\gamma)}+\frac{\partial W}{\partial q^i \partial q^j}\left.\frac{\partial H}{\partial p_j}\right\vert_{im (\gamma)}\right ) dq^i+ \left. \frac{\partial H}{\partial z}\right\vert_{im (\gamma)}dz .
     \end{split}
    \end{equation}
    The right hand side of \eqref{HJT-Con-1-eq-2--} is 
          \begin{equation}   \label{Condu-2} 
  \gamma^*( \mathcal{R}(H)\eta_{\mathcal{Q}}  )   =       \gamma^*\left(\frac{\partial H}{\partial z}(dz- p_i dq^i)\right)= \left. 
          \frac{\partial H}{\partial z}\right\vert_{im (\gamma)}dz-\left. \frac{\partial H}{\partial z}\right\vert_{im (\gamma)}\frac{\partial W}{\partial q^i}dq^i, 
    \end{equation}
    where we have made use of $\eta_{\mathcal{Q}}=dz-p_idq^i$ in the first equality. To retrieve condition \eqref{HJT-Con-1-eq-2--}, we compare the 
    left hand side \eqref{Condu-1} and right hand side \eqref{Condu-2} of \eqref{HJT-Con-1-eq-2--}. The components of $dz$ are the same, therefore they cancel each other. The components of $dq^i$ retrieve precisely the equality \eqref{evo-con-HJ-expli-}. This completes the proof. 
    
    \begin{remark} 
We refer to \cite{de2021hamilton,LeSa17,GrPa20} for a more general picture of the Hamilton-Jacobi theorem for contact Hamiltonian dynamics. See that in the general case the momentum components of the sections $\gamma$ in \eqref{gamma} depend on the fiber coordinate $z$ as well, i.e, $\hat{\gamma}=\hat{\gamma}(q,z)$.
\end{remark}
    
\textbf{Lift of Solutions.} The equations \eqref{evo-con-HJ-expli-} will be referred to as the
 contact  Hamilton-Jacobi equations. Once a solution $W$ is found for \eqref{evo-con-HJ-expli-}, one can lift a solution of the dynamics 
    \begin{equation} \label{HJ-con-red-1}
\frac{d q^i}{dt}=\left. \frac{\partial H}{\partial p_i}\right\vert_{im (\gamma)},\qquad \frac{d z}{dt}= \frac{\partial W}{\partial q^i}\left. \frac{\partial H}{\partial p_i}\right\vert_{im (\gamma)} - \left. H\right\vert_{im (\gamma)}
  \end{equation}
generated by the vector field $X^{\gamma}_H$ to a solution of the dynamics generated by the contact Hamiltonian vector field $X_H$. To see this more explicitly, consider first the dynamics in \eqref{HJ-con-red-1} given by a system of $n+1$ first order differential equations. If $(\varphi^i_t,z_t)$ is a solution of \eqref{HJ-con-red-1}, 
then $(\varphi_t^i,{\partial W}/{\partial q^{i}},z_t)$ is a solution of the contact Hamilton equation \eqref{conham}. Let us see this. The first and the third equation in  \eqref{conham} and the reduced dynamics \eqref{HJ-con-red-1} are the same. So, there only remains to check the second set of equations in \eqref{conham}. A direct calculation shows that the second set $p_i={\partial W \circ \varphi _t}/{\partial q^{i}}$ satisfies
 \begin{equation} \label{solution-2-evo-}
    \begin{split}
 \frac{dp_i}{dt}&=
 \frac{d}{dt}\left(\frac{\partial (W\circ \varphi _t)}{\partial q^{i}}\right)=\frac{\partial^2 (W\circ \varphi _t)}{\partial q^{i}\partial q^{j}}\dot{\varphi}^j_t=\frac{\partial^2 (W\circ \varphi _t)}{\partial q^{i}\partial q^{j}}\frac{\partial H}{\partial p_i}\left(q^i,\frac{\partial (W\circ \varphi _t)}{\partial q^i},z_t\right)
 \\&\qquad =-\frac{\partial H}{\partial q^i}\left(q^i,\frac{\partial (W\circ \varphi _t)}{\partial q^i},z_t\right)
  -\frac{\partial (W\circ \varphi _t)}{\partial q^i}\frac{\partial H}{\partial z}
  \left(q^i,\frac{\partial (W\circ \varphi _t)}{\partial q^i},z_t\right)
 \\&\qquad=-\left.\frac{\partial H}{\partial q^i}\right\vert_{im (\gamma)}-p_i\left.\frac{\partial H}{\partial z}\right\vert_{im (\gamma )}
    \end{split}
 \end{equation}
    where we have employed \eqref{evo-con-HJ-expli-}.  At a glance, \eqref{solution-2-evo-} is precisely the second equation in the contact Hamilton equations. (\ref{conham}). 
    
\subsection{Contact Hamiltonian Dynamics as a Legendrian Submanifold}\label{Sec-Ham-as-Leg}

In this subsection, our goal is to establish a relation between the images of contact Hamiltonian vector fields defined on $\mathcal{T}^*\mathcal{Q}$ and Legendrian submanifolds of the extended tangent bundle $\mathcal{T}\mathcal{T}^*\mathcal{Q}$. For this, we first construct a contact structure on $\mathcal{T}\mathcal{T}^*\mathcal{Q}$ by lifting the contact structure on $\mathcal{T}^*\mathcal{Q}$, see \cite{IbLeMaMa97}.
\begin{theorem}\label{the-tan-con}
For the contact manifold $(\mathcal{T}^*\mathcal{Q},\eta_{\mathcal{Q}})$, the extended tangent bundle $\mathcal{T}\mathcal{T}^*\mathcal{Q} \simeq T\mathcal{T}^*\mathcal{Q}  \times \mathbb{R}$ is a contact manifold with a contact one-form
\begin{equation}\label{eta-T-}
\eta^\mathcal{T}_{\mathcal{Q}} := u  \eta_{\mathcal{Q}} ^V+ \eta_{\mathcal{Q}} ^C
\end{equation}
where $u$ is coordinate on $\mathbb{R}$, whereas $ \eta_{\mathcal{Q}} ^C$ and $ \eta_{\mathcal{Q}} ^V$ are the complete and vertical lifts of $\eta_{\mathcal{Q}}$, respectively. 
\end{theorem} 
Consider the Darboux coordinates $(q^i,p_i,z)$ on $ \mathcal{T}^*\mathcal{Q} $ and the induced coordinates on $\mathcal{T}\mathcal{T}^*\mathcal{Q}\simeq T\mathcal{T}^*\mathcal{Q}\times \mathbb{R}$ as $(q^i,p_i,z,\dot{q}^i,\dot{p}_i,\dot{z},u)$. In this local realization, the lifted contact one-form $ \eta^\mathcal{T}_{\mathcal{Q}}$ defined in \eqref{eta-T-} is computed to be
\begin{equation}\label{eta-T}
\begin{split}
 \eta^\mathcal{T}_{\mathcal{Q}} &=  u  \eta_{\mathcal{Q}} ^V + \eta_{\mathcal{Q}}^C 
  =          u (\dd z - p_i \dd q^i) + (\dd \dot{z} - \dot{p}_i\dd q^i - p_i \dd \dot{q}^i)
  \\ 
  &=\dd \dot{z} +u dz -(\dot{p}_i+u p_i) \dd q^i -  p_i \dd \dot{q}^i,
  \end{split}
        \end{equation}
and the Reeb vector field is $\mathcal{R}^{\mathcal{T}}=\partial/\partial \dot{z}$. 
The one-form $\eta_{\mathcal{Q}}^\mathcal{T}$ is said to be the tangent contact structure.  The image of the section 
\begin{equation}\label{bar-x-H}
(X_H, \mathcal{R}(H)):\mathcal{T}^*\mathcal{Q} \longrightarrow  \mathcal{T}\mathcal{T}^*\mathcal{Q}=T\mathcal{T}^*\mathcal{Q}\times \mathbb{R},  \qquad m\mapsto (X_H(m),\mathcal{R}(H)(m))
\end{equation}
is a Legendrian submanifold of the extended tangent bundle $\mathcal{T}\mathcal{T}^*\mathcal{Q}$. Here, the conformal factor $\mathcal{R}(H)$ is the directional derivative of the Hamiltonian function along the Reeb field.
The inverse of this assertion is also true, that is, if the image space of a section $(X,\lambda)$ of the fibration $\mathcal{T}\mathcal{T}^*\mathcal{Q}\mapsto \mathcal{T}^*\mathcal{Q}$ 
is a Legendrian submanifold of $\mathcal{T}\mathcal{T}^*\mathcal{Q}$, then, at least locally, $X$ is contact Hamiltonian $X=X_H$ and $\lambda=\mathcal{R}(H)$, see \cite{IbLeMaMa97}. For both local and global Hamiltonian sections, their images are  horizontal Legendrian submanifolds. 
We note here that the dynamics determined by a horizontal Legendrian submanifold is given by explicit differential equations. In the present work, our interest is focused in implicit dynamics. For it, we examine non-horizontal  Legendrian submanifolds of $\mathcal{T}\mathcal{T}^*\mathcal{Q}$ in the next subsection.

\subsection{Implicit Contact Hamiltonian Dynamics} \label{impconHam-Sec}

For contact Hamiltonian dynamics on $\mathcal{T}^*\mathcal{Q}$, the image space of the pair $(X_H,\mathcal{R}(H))$ is a horizontal Legendrian submanifold of the tangent contact space $(\mathcal{T}\mathcal{T}^*\mathcal{Q},\eta^\mathcal{T})$. We define an implicit contact Hamiltonian system as a non-horizontal Legendrian submanifold of $\mathcal{T}\mathcal{T}^*\mathcal{Q}$. So that, in order to examine such dynamical systems, we need to explore the geometric foundations of the Legendrian submanifolds of $\mathcal{T}\mathcal{T}^*\mathcal{Q}$. 

As discussed in Subsection \ref{Gen-Fam-Leg}, there exists no section of $\mathcal{T}\mathcal{T}^*\mathcal{Q}\mapsto \mathcal{T}^*\mathcal{Q}$ generating a non-horizontal Legendrian submanifold. Instead, according to Theorem \ref{thm-gen-poin-con}, we can employ a Morse family that plays the role of a generating function. As manifested in Subsection \ref{Gen-Fam-Leg}, a Morse family generates a Legendrian submanifold of the extended cotangent bundle $\mathcal{T}^*\mathcal{T}^*\mathcal{Q}$. To carry this Legendrian submanifold to the tangent contact manifold $\mathcal{T}\mathcal{T}^*\mathcal{Q}$, we need a Morse family with a special contact space \eqref{scs} as in \eqref{Morse-Gen--}. Therefore, to have implicit Hamiltonian dynamics on contact manifolds we need to present a special contact space for $\mathcal{T}\mathcal{T}^*\mathcal{Q}$

\textbf{A Special Contact Space for  $\mathcal{T}\mathcal{T}^*\mathcal{Q}$.} Being an extended cotangent bundle, $\mathcal{T}^*\mathcal{T}^*\mathcal{Q}$ is a contact manifold with canonical contact one-form $\eta_{\mathcal{T}^*\mathcal{Q}}$. We introduce a fiber preserving contact diffeomorphism $\beta^c$ from the tangent contact space $(\mathcal{T}\mathcal{T}^*\mathcal{Q},\eta^{\mathcal{T}}_\mathcal{Q})$ to the extended cotangent bundle $\mathcal{T}\mathcal{T}^*\mathcal{Q}$ as follows
\begin{equation}\label{beta-c}
\beta^c:\mathcal{T}\mathcal{T}^*\mathcal{Q} \longrightarrow \mathcal{T}^*\mathcal{T}^*\mathcal{Q},\qquad (V,u) \mapsto \big (-\iota _ V d \eta _\mathcal{Q} - u \eta_\mathcal{Q}, \eta_\mathcal{Q}(V)\big)
\end{equation}
where $\eta_\mathcal{Q}$ is the canonical contact one-form on $\mathcal{T}^*\mathcal{Q}$. Note that, the following identification 
\begin{equation}\label{contact-Ham-Leg}
\beta^c\circ (X_H,\mathcal{R}(H)) = -\mathcal{T}^* H=-(dH,H)
\end{equation}
holds. Here, $\mathcal{T}^* H$ is the first prolongation of the function $H$. This means that $\beta^c$ maps the horizontal Legendrian submanifold $im(-\mathcal{T}^* H)$ of $\mathcal{T}^*\mathcal{T}^*\mathcal{Q}$ to the horizontal Legendrian submanifold $im(X_H,\mathcal{R}(H))$ of $\mathcal{T}\mathcal{T}^*\mathcal{Q}$. 
Following \eqref{scs-5}, we introduce a special contact space as
\begin{equation}\label{scs-1}
(\mathcal{T}\mathcal{T}^*\mathcal{Q},\tau^0_{\mathcal{T}^*\mathcal{Q}},\mathcal{T}^*\mathcal{Q},\eta^\mathcal{T}_{\mathcal{Q}},\beta^c).
\end{equation}
Here, $\tau^0_{\mathcal{T}^*\mathcal{Q}}:\mathcal{T}\mathcal{T}^*\mathcal{Q}\mapsto \mathcal{T}^*\mathcal{Q}$ is the source projection for the extended tangent bundle.
We merge now this special contact space with a Morse family $-E$ defined on the total space of a fibration $\tau:\mathcal{P}\mapsto \mathcal{T}^*\mathcal{Q}$. So, we plot the following diagram
\begin{equation} \label{Morse-Gen--Con}
\xymatrix{
	\mathbb{R}& \mathcal{P} \ar[dd]^{\tau}\ar[l]^{-E}& \mathcal{T}^*{\mathcal{T}^*\mathcal{Q}} \ar@(ul,ur)^{ \mathcal{N}} \ar[ddr]_{\pi_{\mathcal{T}^*\mathcal{Q}}^0}& & \mathcal{T}\mathcal{T}^*\mathcal{Q}\ar@(ul,ur)^{ \mathcal{N}_{E}} \ar[ll]_{\beta^c} \ar[ddl]^{\tau_{\mathcal{T}^*\mathcal{Q}}^0}\\ \\ &
	\mathcal{T}^*\mathcal{Q} \ar@{=}[rr]&& \mathcal{T}^*\mathcal{Q}
}
\end{equation} 
		where we have employed the cotangent bundle source projection $\pi^0_{\mathcal{T}^*\mathcal{Q}}:\mathcal{T}^*\mathcal{T}^*\mathcal{Q}\mapsto \mathcal{T}^*\mathcal{Q}$.   

Consider the Darboux coordinates $(q^i,p_i,z)$ on $ \mathcal{T}^*\mathcal{Q} $ and assume the induced coordinates on $\mathcal{T}\mathcal{T}^*\mathcal{Q}\simeq T\mathcal{T}^*\mathcal{Q}\times \mathbb{R}$ as $(q^i,p_i,z,\dot{q}^i,\dot{p}_i,\dot{z},u)$. The contact mapping $\beta^c$ in \eqref{beta-c} turns out to be
\begin{equation}\label{beta-c-}
\beta^c(q^i,p_i,z,\dot{q}^i,\dot{p}_i,\dot{z}, u) =(q^i,p_i,z, u p_i+\dot{p}_i ,- \dot{q}^i, -u,\dot{z}-p_i \dot{q}^i).
\end{equation} 
Referring to the commutative diagram in \eqref{Morse-Gen--Con}, a Morse family $-E=-E(q^i,p_i,z,\lambda^a)$, where $\lambda^a$ are the fiber coordinates on $\mathcal{P}$, generates a Legendrian submanifold of $\mathcal{T}^*\mathcal{T}^*\mathcal{Q}$ given by
\begin{equation}\label{Ham-Leg-Sub}
\mathcal{N}=\left\{\left(q^i,p_i,z, -\frac{\partial E}{\partial q^i},-\frac{\partial E}{\partial p_i},-\frac{\partial E}{\partial z},-E \right)\in \mathcal{T}^*\mathcal{T}^*\mathcal{Q}: \frac{\partial E}{\partial \lambda^a} =0 \right\}. 
\end{equation}
Using  the inverse contact diffeomorphism  $(\beta^c)^{-1}$,  we transfer this Legendrian submanifold to a Legendrian submanifold of $\mathcal{T}\mathcal{T}^*\mathcal{Q}$ as follows
\begin{equation}\label{N_E}  
\mathcal{N}_{E} = (\beta^c)^{-1}\left(\mathcal{N}\right) = \left\{\left(q^i,p_i,z, \frac{\partial E}{\partial p_i}, -\frac{\partial E}{\partial z}p_i-\frac{\partial E}{\partial q^i},
p_i\frac{\partial E}{\partial p_i}-E,\frac{\partial E}{\partial z}\right)\in \mathcal{T}\mathcal{T}^*\mathcal{Q}:\frac{\partial E}{\partial \lambda^a} =0 \right\}.  
 \end{equation}
Evidently, this Legendrian submanifold is non-horizontal, hence it describes implicit contact dynamics. In this case, the implicit contact Hamilton equations  generated by a Morse family $E=E(q^i,p_i,z,\lambda^a)$ are just the system of implicit differential equations
\begin{equation}\label{impconham}
\frac{dq^i}{dt}= \frac{\partial E}{\partial p_i}, \qquad \frac{d p _i}{dt} = -\frac{\partial E}{\partial q^i}- 
p_i\frac{\partial E}{\partial z}, \quad \frac{dz}{dt} = p_i\frac{\partial E}{\partial p_i} - E, \qquad \frac{\partial E}{\partial \lambda^a} =0.
\end{equation}

\subsection{HJ Theory for Implicit Contact Hamiltonian Dynamics}\label{Sub-HJ-imp-con}

Theorem \ref{HJT-Con} establishes the Hamilton-Jacobi theory for the explicit Hamiltonian dynamics on the contact framework. In this section, we generalize this result to implicit contact Hamiltonian dynamics \eqref{impconham}. Let $\mathcal{N}_{E}$ be a non-horizontal Legendrian submanifold, generated by a Morse family $E$, of the tangent contact manifold $\mathcal{T}\mathcal{T}^*\mathcal{Q}$. We plot the following commutative diagram
       \begin{equation}\label{TT-evo-}
\xymatrix{
	   &                                                                      & \mathcal{T}\mathcal{T}^*\mathcal{Q} \ar@(ul,ur)^{ \mathcal{N}_{E}} \ar[ldd]_{\psi} \ar[rdd]^{ \tau^0_{\mathcal{T}^*\mathcal{Q}}}&&                                                                     \\  \\
											  & T(\mathcal{Q}   \times \mathbb{R}) \ar[rdd]_{pr_\mathcal{Q}\circ \tau_{\mathcal{Q}   \times \mathbb{R}}}  && \mathcal{T}^*\mathcal{Q} \ar[ldd]^{\pi^0_{\mathcal{Q}}}
										\\	  \\&&\mathcal{Q}
										}
	\end{equation}
	where the projections are given by
	 \begin{equation}
\begin{split}
\psi &: \mathcal{T}\mathcal{T}^*\mathcal{Q}\longrightarrow 
T(\mathcal{Q}   \times \mathbb{R}),\qquad (W,u)\mapsto T\rho(W)\\
pr_\mathcal{Q}\circ \tau_{\mathcal{Q}   \times \mathbb{R}}
&:
T(\mathcal{Q}   \times \mathbb{R})\longrightarrow  \mathcal{Q},\qquad (v,z,\dot{z})\mapsto \tau_{\mathcal{Q}}(v) 
\end{split}	 
	 	\end{equation}
	where $\rho$ is the mapping in \eqref{rho} and $v\in T\mathcal{Q}$. Referring to the commutative diagram, we project the Legendrian submanifold $\mathcal{N}_{E}$  of $\mathcal{T}\mathcal{T}^*\mathcal{Q}$ to a submanifold $C= \tau^0_{\mathcal{T}^*\mathcal{Q}}(\mathcal{N}_{E})$ of the extended cotangent bundle $\mathcal{T}^*\mathcal{Q}$. If the Legendrian submanifold is horizontal then its projection to  $\mathcal{T}^*\mathcal{Q}$ covers the whole manifold. But for a non-horizontal submanifold the projection is only a submanifold of $\mathcal{T}^*\mathcal{Q}$. Next, we project  the Legendrian submanifold $\mathcal{N}_{E}$ to a submanifold $\psi(\mathcal{N}_{E})$ of $T(\mathcal{Q}   \times \mathbb{R})$. This latter projection determines an implicit differential equation on the extended manifold $\mathcal{Q}   \times \mathbb{R}$.  Hamilton-Jacobi theory consists in to retrieve solutions of $\mathcal{N}_{E}$, provided the solutions of the projected submanifold $\psi(\mathcal{N}_{E})$. 
In order to lift the solutions in $\mathcal{Q}   \times \mathbb{R}$ to $\mathcal{T}^*\mathcal{Q}$, we are still in need of a section $\gamma$ as given in \eqref{gamma}, but two ingredients
of the theory are missing. One is that the base manifold, denoted by $C= \tau^0_{\mathcal{T}^*\mathcal{Q}}(\mathcal{N}_{E})$, is not necessarily the whole $\mathcal{T}^*Q$, but possibly a proper submanifold of it.
The second is the nonexistence of a Hamiltonian vector field due to the
implicit character of the equations. 

Let $\gamma$ be a section of the fibration $\mathcal{T}^*\mathcal{Q}\mapsto \mathcal{Q}   \times \mathbb{R}$ in the form of \eqref{gamma} and then consider the restriction of the Legendrian submanifold $\mathcal{N}_{E}$ to the image of $\gamma$. This reads the following restricted submanifold 
\begin{equation} \label{N-Im-E}
\begin{split}
\mathcal{N}
_{E}\vert_{im(\gamma)}&=\Big \{\Big(q^i,\frac{\partial W}{\partial q^i},z; \left.\frac{\partial E}{\partial p_i}\right\vert_{im(\gamma)},\left.
-\frac{\partial E}{\partial q^i}\right\vert_{im(\gamma)}-\left.\frac{\partial W}{\partial q^i}\frac{\partial E}{\partial z}\right\vert_{im(\gamma)},\\ &\hspace{2.5cm}
\frac{\partial W}{\partial q^i}
\left.
\frac{\partial E}{\partial p_i}\right\vert_{im(\gamma)}-\left. E \right\vert_{im(\gamma)},  \left.\frac{\partial E}{\partial z}\right\vert_{im(\gamma)}\Big )
\in \mathcal{T}\mathcal{T}^*\mathcal{Q}  :\left.\frac{\partial E}{\partial \lambda^a}\right\vert_{im(\gamma)}=0\Big\}.
\end{split}
\end{equation}
If we are in the explicit case, that is if the submanifold $\mathcal{N}
_{E}$ is the image of the two-tuple $(X_H,\mathcal{R}(H))$, then $\mathcal{N}
_{E}\vert_{im(\gamma)}$ reduces to the image space of $(X_H,\mathcal{R}(H)) \circ \gamma$. The submanifold $ \mathcal{N}
_{E}\vert_{im(\gamma)}$ depends only on the coordinates $(q^i,z)$. So that, by the mapping $\psi$, we can project it to  a submanifold $\mathcal{N}_E^{\gamma}$ of $T(\mathcal{Q}\times \mathbb{R})$ as
\begin{equation} \label{N-gamma}
\begin{split}
\mathcal{N}_E^{\gamma}=\psi (\mathcal{N}
_{E}\vert_{im(\gamma)})&=\Big \{\Big(q^i,z;\left.\frac{\partial E}{\partial p_i}\right\vert_{im(\gamma)},
\frac{\partial W}{\partial q^i}
\left.
\frac{\partial E}{\partial p_i}\right\vert_{im(\gamma)}-\left. E \right\vert_{im(\gamma)} 
\Big)\\&\hspace{2.5cm} \in T(\mathcal{Q}\times \mathbb{R}):\left.\frac{\partial E}{\partial \lambda^a}\right\vert_{im(\gamma)}=0 \Big\}.
\end{split}
\end{equation}
In general, the submanifold $\mathcal{N}_E^{\gamma}$ defines an implicit differential equation on $\mathcal{Q}\times \mathbb{R}$ given as
\begin{equation}
\frac{dq^i}{dt}	=\left.\frac{\partial E}{\partial p_i}\right\vert_{im(\gamma)},\qquad 
\frac{dz}{dt}	 =
\frac{\partial W}{\partial q^i}
\left.
\frac{\partial E}{\partial p_i}\right\vert_{im(\gamma)}-\left. E \right\vert_{im(\gamma)} , \left.\qquad \frac{\partial E}{\partial \lambda^a}\right \vert_{im(\gamma)}=0. 
 \end{equation}
 Now, we are ready to state Hamilton-Jacobi theorem for the implicit contact Hamiltonian dynamics generalizing Theorem \ref{HJT-Con}. 
   
\begin{theorem} \label{HJ-imp-con}
For a section $\gamma$ in the form \eqref{gamma} such that $\hat{\gamma}$ is closed, the following two conditions are equivalent:
\begin{enumerate}
\item The submanifolds $\mathcal{N}
_{E}$ in \eqref{N_E} and $\mathcal{N}
_{E}^{\gamma}$ in \eqref{N-gamma} are '$\gamma$-related', that is
\begin{equation}\label{HJT-Con-1-eq-1-}
\big( T\gamma (\mathcal{N}
_{E}^{\gamma}), \mathcal{R}(E)\big) = \mathcal{N}
_{E}\vert_{im(\gamma)}
\end{equation}
where $T\gamma:T(\mathcal{Q}\times \mathbb{R})\mapsto T\mathcal{T}^{*}\mathcal{Q}$ is tangent mapping of  the section $\gamma$, and $ \mathcal{R}$ is the Reeb field. $ \mathcal{R}(E)$ is the directional derivative $\partial E / \partial z$ of the Morse family in the direction of the component $z$.
\item The equation
\begin{equation}\label{HJT-Con-1-eq-2-}
d( E \circ \gamma)=  \gamma^*(\mathcal{R}(E)\eta_{\mathcal{Q}})
\end{equation}
 is fulfilled. Here, $E \circ \gamma(q,z,\lambda)=E(q^i,\partial W / \partial q^i, z,\lambda^a)$.
\end{enumerate}
\end{theorem} 
  In local coordinates, the implicit Hamilton-Jacobi equations for contact Hamiltonian dynamics is computed to be
  \begin{equation}\label{HJ-imp-exp-cont}
  \frac{\partial W}{\partial q^i \partial q^j}\left. \frac{\partial E}{\partial p_j}\right\vert_{im(\gamma)} +\left. \frac{\partial E}{\partial q^i}\right\vert_{im(\gamma)}
  +\frac{\partial W}{\partial q^i}\left. \frac{\partial E}{\partial z}
  \right\vert_{im(\gamma)}=0, \qquad 
\left.  \frac{\partial E}{\partial \lambda^a}\right \vert_{im(\gamma)}=0. 
  \end{equation}
  \section{Implicit Evolution Hamiltonian Dynamics and Hamilton--Jacobi Theory}\label{Sec-Imp-Evo}

\subsection{Evolution Contact Dynamics as a Lagrangian Submanifold}\label{Sec-Ev-asLag}

Consider the contact manifold $(\mathcal{T}^{*}\mathcal{Q},\eta_{\mathcal{Q}})$.  We cite \cite{esen2021contact} for  all of the results in this Subsection. In \eqref{decomp-TM}, we have stated that the kernel of the contact one-form is a symplectic vector bundle of the tangent bundle. Define the inclusion mapping
	\begin{equation}\label{inc}
	j = (i,{\mathrm{Id}}_{\mathbb{R}}):\mathfrak{H} \mathcal{T}^{*}\mathcal{Q} \times \mathbb{R} \hookrightarrow \mathcal{T}\mathcal{T}^{*}\mathcal{Q} = T\mathcal{T}^{*}\mathcal{Q} \times \mathbb{R}
	\end{equation}
	where $i$ is the inclusion of $\mathfrak{H} \mathcal{T}^{*}\mathcal{Q} $ into $T\mathcal{T}^{*}\mathcal{Q} $. We pull-back the tangent contact one-form $\eta^\mathcal{T}_{\mathcal{Q}}$ in \eqref{eta-T-} to $\mathfrak{H} \mathcal{T}^{*}\mathcal{Q} \times \mathbb{R}$ by the mapping $j$. So, we obtain a one-form 
$\theta_{\eta} = j^* \eta^\mathcal{T}_{\mathcal{Q}}$ on $\mathfrak{H} \mathcal{T}^{*}\mathcal{Q} \times \mathbb{R}$.
\begin{theorem}\label{the-symp}
The exact two-form $\omega_{\eta} =   d \theta_{\eta}$ induces a symplectic structure on $ \mathfrak{H} \mathcal{T}^{*}\mathcal{Q} \times \mathbb{R}$. 
\end{theorem}
{\rm\textbf{Proof.}} Let $(q^i,p_i,z)$ be the Darboux coordinates on $\mathcal{T}^{*}\mathcal{Q}$, so that  $(q^i,p_i,z, \dot{q}^i,\dot{p}_i,\dot{z}, u)$ are the induced  coordinates for the iterated extended  tangent bundle $\mathcal{T}\mathcal{T}^{*}\mathcal{Q}$. We employ $(q^i,p_i,z, \dot{q}^i,\dot{p}_i, u)$ as a local coordinate chart on $\mathfrak{H} \mathcal{T}^{*}\mathcal{Q} \times \mathbb{R}$ since
	\begin{equation}\label{map-j}
		j(q^i,p_i,z, \dot{q}^i,\dot{p}_i, u) = (q^i,p_i,z, \dot{q}^i,\dot{p}_i, p_i \dot{q}^i, u).
	\end{equation}
The Reeb vector field $\mathcal{R}^{\mathcal{T}}=\partial/\partial \dot{z}$ of the contact manifold $(\mathcal{T}\mathcal{T}^{*}\mathcal{Q},\eta^\mathcal{T}_{\mathcal{Q}})$ is transverse to the manifold of codimension one $\mathfrak{H} \mathcal{T}^{*}\mathcal{Q} \times \mathbb{R}$. This implies that the exact two-form $\omega_{\eta} = d \theta_{\eta}=d(j^* \eta^\mathcal{T}_{\mathcal{Q}})$ induces a symplectic structure on  $\mathfrak{H} \mathcal{T}^{*}\mathcal{Q} \times \mathbb{R}$. Locally, one has that
	\begin{equation}\label{omega-H}
	\begin{split}
		\theta_{\eta} &= 
		u dz -(\dot{p}_i+u p_i) \dd q^i + \dot{q}^i \dd  p_i,
\\
		\omega_{\eta}& = du \wedge d z - d \dot{p}_i \wedge \dd q^i - p_i du \wedge \dd q^i - u d p_i \wedge \dd q^i  + d \dot{q}^i \wedge d  p_i.
\end{split}
	\end{equation}

	\textbf{Lagrangian submanifolds and evolution vector fields.} Given a vector field $X$ on $\mathcal{T}^{*}\mathcal{Q}$ and a real smooth function $f:\mathcal{T}^{*}\mathcal{Q}\to \mathbb{R}$, one can construct a section $(X,f): \mathcal{T}^{*}\mathcal{Q} \to \mathcal{T} \mathcal{T}^{*}\mathcal{Q}$ of the extended tangent bundle $\mathcal{T} \mathcal{T}^{*}\mathcal{Q}$.
	\begin{theorem}\label{lag-evo-char}
	Let $(\mathcal{T}^{*}\mathcal{Q},\eta_{\mathcal{Q}})$ be a contact manifold. 
	The map	$(X,f):\mathcal{T}^{*}\mathcal{Q}\mapsto \mathcal{T}\mathcal{T}^{*}\mathcal{Q}$ defines a Lagrangian submanifold of the exact symplectic manifold $(\mathfrak{H} \mathcal{T}^{*}\mathcal{Q} \times \mathbb{R}, \omega_{\eta})$ if and only if $\eta_{\mathcal{Q}}(X) = 0$ and $\mathcal{L}_X \eta_{\mathcal{Q}} + f \eta_{\mathcal{Q}}$ is closed.
	\end{theorem}

	{\rm\textbf{Proof.}   }
	First of all, note that the image of $X$ lies on $\mathfrak{H} \mathcal{T}^{*}\mathcal{Q} \times \mathbb{R}$ if and only if $\eta_{\mathcal{Q}}(X) = 0$.
	For the second condition we use well-known properties of complete and vertical lifts (see \cite{YaIs73});
	\begin{equation}
		(X, f)^* \eta_{\mathcal{Q}}^\mathcal{T} = 
		(X, f)^* (u \eta_{\mathcal{Q}}^V + \eta_{\mathcal{Q}}^C) = 
		f X^* \eta_{\mathcal{Q}}^V + X^* \eta^C = f \eta _{\mathcal{Q}}+ \mathcal{L}_X \eta_{\mathcal{Q}}.
	\end{equation}
	Thus, the image of $(X, f)$ is Lagrangian if and only if
	\begin{equation}
		(X, f)^*\omega_{\eta} = d (X, f)^* \eta^\mathcal{T}_{\mathcal{Q}} =
		d( f \eta_{\mathcal{Q}} + \mathcal{L}_X (\eta_{\mathcal{Q}})) = 0.
	\end{equation}

	Looking at the definition of the evolution vector field~\eqref{evo-def}, we obtain the following result.
	\begin{corollary}\label{corola}
		The map	$(X,f):\mathcal{T}^{*}\mathcal{Q}\mapsto \mathcal{T}\mathcal{T}^{*}\mathcal{Q}$ defines a Lagrangian submanifold of the exact symplectic manifold $(\mathfrak{H} \mathcal{T}^{*}\mathcal{Q} \times \mathbb{R}, \omega_{\eta})$ if and only if, locally, $X = \varepsilon_H$ and $f = \mathcal{R}(H)$ for a (local) smooth function $H: \mathcal{T}^{*}\mathcal{Q} \to \mathcal{R}$.
	\end{corollary}
	\textbf{Proof.}
	Since the evolution vector field satisfies
	\begin{equation}
		\mathcal{L}_{\varepsilon_H}\eta_{\mathcal{Q}}=dH-\mathcal{R}(H)\eta_{\mathcal{Q}},\qquad \eta_{\mathcal{Q}}(\varepsilon_H)=0,
	\end{equation}
one deduces that $\mathcal{L}_X \eta_{\mathcal{Q}} + f \eta_{\mathcal{Q}} = d H$ is closed.
	Conversely, if $\eta(X) = 0$ and $\mathcal{L}_X \eta_{\mathcal{Q}} + f \eta_{\mathcal{Q}}$ is closed, in a local chart $U$, one has that  
	\begin{equation}
		\left.\mathcal{L}_X \eta _{\mathcal{Q}}+ f \eta_{\mathcal{Q}}\right \vert_U = \left. \contr{X} d \eta_{\mathcal{Q}} + f \eta_{\mathcal{Q}}\right \vert_U  = d H 
	\end{equation}
	for a local smooth  function $H$. By contracting the equality with the Reeb vector field, we obtain that $f = \mathcal{R}(H)$. Thus, $X = \varepsilon_H$.
	
We consider now the following commutative diagram:
\begin{equation}
	\begin{tikzcd}
		{(\mathcal{T}\mathcal{T}^{*}\mathcal{Q},\eta_{\mathcal{Q}}^{\mathcal{T}})} \arrow[rr, "\beta^c"]                      &                                                                      & {(\mathcal{T}^*\mathcal{T}^{*}\mathcal{Q}, \eta_{\mathcal{T}^{*}\mathcal{Q}})}            \\ \\
		{(\mathfrak{H} \mathcal{T}^{*}\mathcal{Q}  \times \mathbb{R}, \omega_{\eta})} \arrow[uu, "j", hook] \arrow[rr, "\beta^0"] &                                                                      & {(T^*\mathcal{T}^{*}\mathcal{Q}, \omega_{\mathcal{T}^{*}\mathcal{Q}})}. \arrow[uu, "{(\mathrm{Id}_{T\mathcal{T}^{*}\mathcal{Q}},0)}",swap, hook] 		\end{tikzcd}
	\end{equation}
	Here,  $\beta^c$ is the contact mapping in \eqref{beta-c}, the mapping $(\mathrm{Id}_{T\mathcal{T}^{*}\mathcal{Q}},0)$ is the canonical inclusion of the cotangent bundle into the extended cotangent bundle as a zero section and the smooth map $\beta^0$ is given by
	\begin{equation}\label{beta-0}
	 \beta^0:\mathfrak{H} \mathcal{T}^{*}\mathcal{Q}  \times \mathbb{R}\longrightarrow T^*\mathcal{T}^{*}\mathcal{Q},\qquad (V,u)\mapsto -\contr{V} d \eta_{\mathcal{Q}} - u \eta_{\mathcal{Q}}.
	 \end{equation}
In coordinates, 
\begin{equation}\label{beta-0-loc}
\beta^0(q^i,p_i,z,\dot{q}^i,\dot{p}_i, u) =(q^i,p_i,z, u p_i+\dot{p}_i ,- \dot{q}^i, -u).
	 \end{equation}
We plot the following diagram
\begin{equation}\label{Right-w-evo}
	\begin{tikzcd}
		\mathfrak{H} \mathcal{T}^{*}\mathcal{Q}   \times \mathbb{R} \arrow[rr, "\beta^0"] \arrow[rdd,"\hat{\tau}_{\mathcal{T}^*\mathcal{Q}}"]&                                                                      & T^*\mathcal{T}^*\mathcal{Q} \arrow[ldd,swap,"\pi_{\mathcal{T}^*\mathcal{Q}}"] \\  \\
										  & \mathcal{T}^*\mathcal{Q} \arrow[luu, bend left, "{(\varepsilon_H, \mathcal{R}(H))}"] \arrow[ruu, bend right, swap,"-d H"] &                                                             
		\end{tikzcd}
	\end{equation}
where $\hat{\tau}_{\mathcal{T}^*\mathcal{Q}}$ is the projection taking a pair $(V,u)$ mapping to $\tau_{\mathcal{T}^*\mathcal{Q}}(V)$, that is,
 \begin{equation}
 \hat{\tau}_{\mathcal{T}^*\mathcal{Q}}:\mathfrak{H} \mathcal{T}^{*}\mathcal{Q}  \times \mathbb{R}\longrightarrow \mathcal{T}^*\mathcal{Q},\qquad (V,u)\mapsto \tau_{\mathcal{T}^*\mathcal{Q}}(V)
 \end{equation}
using that $ \mathfrak{H} \mathcal{T}^{*}\mathcal{Q}  $ is a vector subbundle of $ T\mathcal{T}^*\mathcal{Q} $. In order to see that the triangle commutes, we compute that
\begin{equation}
	\beta^0\circ (\varepsilon_H,\mathcal{R}(H)) = - \contr{\varepsilon_H} \dd \eta_\mathcal{Q} - \mathcal{R}(H) \eta_\mathcal{Q} = -\mathcal{L}_{\varepsilon_H} \eta_\mathcal{Q}- \mathcal{R}(H) \eta_\mathcal{Q} = - d H,
\end{equation}
where we have used~\eqref{evo-def} and Cartan's formula.  
Since $\beta^0$ is a symplectic diffeomorphism, the image  $im(\varepsilon_H,\mathcal{R}(H))$ is a Lagrangian submanifold of $\mathfrak{H} \mathcal{T}^{*}\mathcal{Q}   \times \mathbb{R}$. 

\subsection{HJ Theory for Evolutionary Contact Dynamics}\label{Sub-HJ-evol}

We consider a real valued function $W$ on the base manifold $\mathcal{Q}$. Its first prolongation $\mathcal{T}^*W$ is defined in \eqref{j1F}. Recall that the image of $\mathcal{T}^*W$ is a Legendrian submanifold of the extended cotangent space $\mathcal{T}^{*}\mathcal{Q}$. Consider the following diagram
 \begin{equation}
\xymatrix{ \mathcal{T}^{*}\mathcal{Q}\simeq T^{*}\mathcal{Q}\times \mathbb{R}
\ar[dd]^{\pi^0_\mathcal{Q}} \ar[rrr]^{\varepsilon_H}&   & &T\mathcal{T}^{*}Q\ar[dd]^{T\pi^0_\mathcal{Q}}\\
  &  & &\\
\mathcal{Q}  \ar@/^2pc/[uu]^{\mathcal{T}^*W}\ar[rrr]^{\varepsilon^{\mathcal{T}^*W}_H}&  & & T\mathcal{Q}}
\end{equation}
where $\pi^0_\mathcal{Q}$ is the projection defined in \eqref{Sec-Jet}. Referring to this commutative diagram and 
 starting with  an evolution Hamiltonian vector field $\varepsilon_H$ defined on   $\mathcal{T}^{*}Q$, we define a vector field $\varepsilon^{\mathcal{T}^*W}_H$ on $\mathcal{Q}$ as
 \begin{equation}
\varepsilon^{\mathcal{T}^*W}_H= T\pi^0_\mathcal{Q} \circ  X_H \circ \mathcal{T}^*W
\end{equation}
where $T\pi^0_\mathcal{Q}$ is the tangent mapping of $\pi^0_\mathcal{Q}$. The following is a Hamilton-Jacobi theorem for evolution Hamiltonian dynamics.
\begin{theorem} \label{HJT-Evo-Con-1}
For a smooth function $W=W(q)$ on $\mathcal{Q}$, the following conditions are equivalent:
\begin{enumerate}
\item The vector fields $\varepsilon_{H}$ and $\varepsilon^{\mathcal{T}^*W}_H$ are $\mathcal{T}^*W$-related, that is
\begin{equation}\label{HJT-Con-1-eq-1}
T\mathcal{T}^*W \circ  \varepsilon_H^{\mathcal{T}^*W} =\varepsilon_H\circ \mathcal{T}^*W,
\end{equation}
where $T\mathcal{T}^*W:T\mathcal{Q}\mapsto T\mathcal{T}^{*}\mathcal{Q}$ is the tangent mapping of the prolongation $\mathcal{T}^*W$. 
\item The equation 
\begin{equation}\label{HJT-Con-1-eq-2}
d( H\circ \mathcal{T}^*W  )=0,
\end{equation}
is fulfilled.
\end{enumerate}
\end{theorem}

\textbf{Proof.}
Let us start the proof by determining the local realization of the vector field 
 \begin{equation}
 \varepsilon^{\mathcal{T}^*W}_H: \mathcal{Q} \longrightarrow T\mathcal{Q},\qquad (q^i)\mapsto \Big(q^i,\left. \frac{\partial H}{\partial p_i}\right\vert_{im \mathcal{T}^*W} \Big).
   \end{equation}
We now compute the left and the right hand sides of \eqref{HJT-Con-1-eq-1}. In induced coordinates $(q^i,\dot{q}^i)$ on the tangent bundle $T\mathcal{Q}$, the tangent mapping of $\mathcal{T}^*W$ is
 \begin{equation}
 T\mathcal{T}^*W(q^i,\dot{q}^i)=\left(q^i,\frac{\partial W}{\partial q^i},W;\dot{q}^i,\frac{\partial W}{\partial q^i \partial q^j}\dot{q}^j,\frac{\partial W}{\partial q^i}\dot{q}^i\right)
 \end{equation}
then the left hand side of the equation \eqref{HJT-Con-1-eq-1} is computed to be
 \begin{equation}\label{Proof-1}
 \begin{split}
 T\mathcal{T}^*W\circ   \varepsilon_H^{\mathcal{T}^*W} (q^i)&=\Big(q^i,\frac{\partial W}{\partial q^i},W; \left. \frac{\partial H}{\partial p_i}\right\vert_{im \mathcal{T}^*W}, \frac{\partial W}{\partial q^i \partial q^j}\left. \frac{\partial H}{\partial p_j}\right\vert_{im \mathcal{T}^*W},\frac{\partial W}{\partial q^i}\left. \frac{\partial H}{\partial p_i}\right\vert_{im \mathcal{T}^*W}\Big).
 \end{split}
  \end{equation}
On the other hand, the right hand side of \eqref{HJT-Con-1-eq-1}  is
  \begin{equation}\label{Proof-2}
   \begin{split}
 \varepsilon_H\circ \mathcal{T}^*W(q^i)&=\Big(q^i,\frac{\partial W}{\partial q^i},W ;
  \left. \frac{\partial H}{\partial p_i}\right\vert_{im \mathcal{T}^*W}, -\left. \frac{\partial H}{\partial q^i}\right\vert_{im \mathcal{T}^*W}
  -\frac{\partial W}{\partial q^i}\left. \frac{\partial H}{\partial z}
 \right\vert_{im \mathcal{T}^*W},
  \frac{\partial W}{\partial q^i}\left.
  \frac{\partial H}{\partial p_i}\right\vert_{im \mathcal{T}^*W}
  \Big). 
   \end{split}
    \end{equation}
According to the identity  \eqref{HJT-Con-1-eq-1}, by equating \eqref{Proof-1} and \eqref{Proof-2} one arrives at that the basis points, that are the first three entries of the local realizations, are the same. Two of the fiber variables, that are the fourth and sixth entries, are also the same. The fifth entries of \eqref{Proof-1} and \eqref{Proof-2}  are not the same. They are equal if and only if the following identity holds
 \begin{equation}\label{evo-con-HJ-expli}
\frac{\partial W}{\partial q^i \partial q^j}\left. \frac{\partial H}{\partial p_j}\right\vert_{im \mathcal{T}^*W}+\left. \frac{\partial H}{\partial q^i}\right\vert_{im \mathcal{T}^*W}
  +\frac{\partial W}{\partial q^i}\left. \frac{\partial H}{\partial z}
  \right\vert_{im \mathcal{T}^*W}=0
     \end{equation}
We can write this identity as 
\begin{equation}
d  H \left(q^i,\frac{\partial W}{\partial q^i},W\right) =0
\end{equation}
 which is precisely the condition given in \eqref{HJT-Con-1-eq-2}. For the proof of the inverse assertion, one only needs to follow the same steps but in the reversed order. This completes the proof. 
 
In terms of the Darboux' coordinates the Hamilton-Jacobi equation \eqref{HJT-Con-1-eq-2} turns out to be 
 \begin{equation}\label{evo-con-HJ}
H \big(q^i,\frac{\partial W}{\partial q^i},W \big)=c
 \end{equation}
 where $c$ is the constant of integration. We call \eqref{evo-con-HJ} as evolution Hamilton-Jacobi equation.  
 
\textbf{Lifts of Solutions.}
  Once a solution $W$ is found for the evolution  Hamilton-Jacobi equation \eqref{evo-con-HJ}, we can lift the solutions of the projected dynamics $\varepsilon^{\mathcal{T}^*W}_H$ on $\mathcal{Q}$ to the solutions of the evolution dynamics $\varepsilon_H$ on the extended cotangent bundle $T^{*}\mathcal{Q}$ by means of the first prolongation of $W$. To see this, consider first the dynamics $\dot{q}=\varepsilon^{\mathcal{T}^*W}_H$ given in local coordinates as 
 \begin{equation} \label{redHamEq-Int-evo}
\dot{q}^i=\left. \frac{\partial H}{\partial p_i}\right\vert_{im \mathcal{T}^*W}.
\end{equation}
Notice that, the system \eqref{redHamEq-Int-evo} consists of a number $n$ of first order differential equations, whereas (\ref{evo-eq}) has $2n+1$ equations. Accordingly, it is easier to solve this system than solving the Hamilton equations. If $W$ is a solution of the Hamilton-Jacobi problem, then a solution of \eqref{redHamEq-Int-evo} can be lifted to a solution of the evolution Hamilton equations. More concretely, if $\varphi_t=(\varphi^i_t)$ is a solution of \eqref{redHamEq-Int-evo} 
then $(\varphi_t,{\partial W}/{\partial q^{i}}(\varphi _t),W(\varphi_t))$ is a solution of the evolution contact Hamilton equations \eqref{evo-eq}. Indeed a direct calculation shows that the second term $p_i={\partial W}/{\partial q^{i}}(\varphi _t)$ in the proposed solution satisfies
 \begin{equation} \label{solution-2-evo}
    \begin{split}
 \frac{dp_i}{dt}&=
 \frac{d}{dt}\left(\frac{\partial (W\circ \varphi _t)}{\partial q^{i}}\right)=\frac{\partial^2 (W\circ \varphi _t)}{\partial q^{i}\partial q^{j}}\dot{\varphi}^j_t=\frac{\partial^2 (W\circ \varphi _t)}{\partial q^{i}\partial q^{j}}\frac{\partial H}{\partial p_i}\left(q^i,\frac{\partial (W\circ \varphi _t)}{\partial q^i},(W\circ \varphi _t)(q)\right)
 \\&\qquad =-\frac{\partial H}{\partial q^i}\left(q^i,\frac{\partial (W\circ \varphi _t)}{\partial q^i},(W\circ \varphi _t)(q)\right)
  -\frac{\partial (W\circ \varphi _t)}{\partial q^i}\frac{\partial H}{\partial z}
  \left(q^i,\frac{\partial (W\circ \varphi _t)}{\partial q^i},(W\circ \varphi _t)(q)\right)
 \\&\qquad=-\left.\frac{\partial H}{\partial q^i}\right\vert_{\mathcal{T}^*(W\circ \varphi _t)}-p_i\left.\frac{\partial H}{\partial z}\right\vert_{\mathcal{T}^*(W\circ \varphi _t)}
    \end{split}
 \end{equation}
    where in the first line we have employed the solution $\varphi_t=(\varphi^i_t)$ and in the second line, we  have substituted the explicit Hamilton-Jacobi equation in \eqref{evo-con-HJ-expli}. Looking carefully, one concludes that \eqref{solution-2-evo} is precisely the second equation in the evolution contact Hamilton equation (\ref{evo-eq}). We consider $z=W\circ \varphi _t$ and take the time derivative of it that is
     \begin{equation} 
         \begin{split}
    \frac{dz}{dt}& =
 \frac{d}{dt}\big(W\circ \varphi _t\big)=\frac{\partial  (W\circ \varphi _t)}{\partial q^{i}}\dot{\varphi}^i_t =\frac{\partial  (W\circ \varphi _t)}{\partial q^{i}}\frac{\partial H}{\partial p_i}\left(q^i,\frac{\partial (W\circ \varphi _t)}{\partial q^i},(W\circ \varphi _t)(q)\right)
  \\& =p_i\left.\frac{\partial H}{\partial p_i}\right\vert_{\mathcal{T}^*(W\circ \varphi _t)}
  \end{split}
     \end{equation}
   where in the first line we have employed the solution $\varphi_t=(\varphi^i_t)$ and $p_i={\partial W}/{\partial q^{i}}(\varphi _t)$. See that, this is the third equation in the evolution contact Hamilton equation (\ref{evo-eq}). 
    
    \subsection{Implicit Evolution Hamiltonian Dynamics}
\label{Imp-Evo-Sec}
    
    As manifested in the previous section, for a given Hamiltonian function $H$ on the contact space $\mathcal{T}^*\mathcal{Q}$, the image of the pair $(\varepsilon_H,\mathcal{R}(H))$ is a (horizontal) Lagrangian submanifold of the symplectic manifold  $\big(\mathfrak{H} \mathcal{T}^{*}\mathcal{Q}  \times \mathbb{R},\omega_\eta\big)$. Evidently, not all Lagrangian submanifolds of $\mathfrak{H} \mathcal{T}^{*}\mathcal{Q}  \times \mathbb{R}$ are horizontal. We define implicit evolution Hamiltonian dynamics  as a non-horizontal Lagrangian submanifold of $\mathfrak{H} \mathcal{T}^{*}\mathcal{Q}  \times \mathbb{R}$. 
    
    Referring to the special symplectic structure in \eqref{Right-w-evo}, we can identify the Lagrangian submanifolds of  $\mathfrak{H} \mathcal{T}^{*}\mathcal{Q}    \times \mathbb{R}$ and the Lagrangian submanifolds of the cotangent bundle $T^*\mathcal{T}^*\mathcal{Q}$. In this picture, the symplectic diffeomorphism $\beta^0$ in \eqref{beta-0} is fiber preserving so that it maps horizontal Lagrangian submanifolds to  horizontal ones, and non-horizontal Lagrangian submanifolds to  non-horizontal ones. In the light of the Maslov-Hörmander theorem (see Appendix \ref{appendix-1}), we can argue that for any (including non-horizontal) Lagrangian submanifold $T^*\mathcal{T}^*\mathcal{Q}$ there exists a Morse family $-F$ defined on the total space of a fiber bundle $P\mapsto \mathcal{T}^*\mathcal{Q}$. We denote this Lagrangian submanifold by $\mathcal{S}$. The inverse of the mapping $\beta^0$ maps this non-horizontal Lagrangian submanifold $\mathcal{S}_{F}=(\beta^0)^{-1}(\mathcal{S})$ to a  non-horizontal Lagrangian submanifold of $\mathfrak{H} \mathcal{T}^{*}\mathcal{Q}   \times \mathbb{R}$. Since, $\beta^0$  is a diffeomorphism, this holds for all non-horizontal Lagrangian submanifold of $\mathfrak{H} \mathcal{T}^{*}\mathcal{Q}    \times \mathbb{R}$. We depict this in the following diagram:
     \begin{equation}\label{Right-w-evo-merge}
\xymatrix{
		\mathfrak{H} \mathcal{T}^{*}\mathcal{Q}  \times \mathbb{R} \ar@(ul,ur)^{ \mathcal{S}_{F}} \ar[rr]^{\beta^0} \ar[rrdd]^{\hat{\tau}_{\mathcal{T}^*\mathcal{Q}}}&                                                                      & \ar@(ul,ur)^{ \mathcal{S}} T^*\mathcal{T}^*\mathcal{Q} \ar[dd] ^{\pi_{\mathcal{T}^*\mathcal{Q}}"} &&                                                               P \ar[dd]\ar[rr]^{-F}&&                                                              \mathbb{R} \\  \\
									  && \mathcal{T}^*\mathcal{Q} \ar@{=}[rr]&&                                                               \mathcal{T}^*\mathcal{Q}
									  }
	\end{equation}
       In coordinates, assuming the induced coordinates $(q^i,p_i,z,\lambda^a)$ on the total space $P$, the Lagrangian $\mathcal{S}$ is computed to be
   \begin{equation} \label{S-Loc-Evo}
\mathcal{S}=\left \{\left(q^i,p_i,z;-\frac{\partial F}{\partial q^i},-\frac{\partial F}{\partial p_i},-\frac{\partial F}{\partial z}\right )\in T^*\mathcal{T}^*Q:\frac{\partial F}{\partial \lambda^a}=0\right\}
\end{equation}
Referring to the local realization of the mapping $\beta^0$ given in \eqref{beta-0-loc} we compute the Lagrangian submanifold $\mathcal{S}_{F}$ as 
   \begin{equation} \label{S-Loc-Evo-F}
\mathcal{S}_{F}=\left \{\left(q^i,p_i,z;\frac{\partial F}{\partial p_i},-\frac{\partial F}{\partial q^i}-p_i\frac{\partial F}{\partial z},\frac{\partial F}{\partial z}\right )\in \mathfrak{H} \mathcal{T}^{*}\mathcal{Q}  \times \mathbb{R}:\frac{\partial F}{\partial \lambda^a}=0\right\}.
\end{equation}
To write this as a system of implicit differential equations, we map the Lagrangian submanifold $\mathcal{S}_{F}$ into the extended tangent bundle $\mathcal{T}\mathcal{T}^*\mathcal{Q}$ by the mapping $j$ in \eqref{map-j}. We arrive at
   \begin{equation} \label{S-Loc-Evo-F-j}
j(\mathcal{S}_{F})=\left \{\left(q^i,p_i,z;\frac{\partial F}{\partial p_i},-\frac{\partial F}{\partial q^i}-p_i\frac{\partial F}{\partial z},p_i\frac{\partial F}{\partial p_i},\frac{\partial F}{\partial z}\right )\in  \mathcal{T}\mathcal{T}^*\mathcal{Q}:\frac{\partial F}{\partial \lambda^a}=0\right\}.
\end{equation}
Accordingly, the implicit evolution dynamics generated by the Morse family $F=F(q^i,p_i,z,\lambda^a)$ is 
\begin{equation}\label{evo-eq-imp}
\frac{dq^i}{dt}	= \frac{\partial F}{\partial p_i}, \qquad \frac{dp_i}{dt}	= -\frac{\partial F}{\partial q^i}- 
	p_i\frac{\partial F}{\partial z}, \quad \frac{dz}{dt} = p_i\frac{\partial F}{\partial p_i}, \qquad \frac{\partial F}{\partial \lambda^a}=0.
\end{equation}

\subsection{HJ Theory for Implicit Evolution Hamiltonian Dynamics}\label{Sub-HJ-imp-evo-11}

The (evolution) Hamilton-Jacobi Theorem \ref{HJT-Evo-Con-1} only deals with explicit evolution Hamiltonian systems. In this section, we present a generalization of this theorem for implicit evolution Hamiltonian dynamics introduced in Subsection \ref{Imp-Evo-Sec}. 

Consider the following diagram     
    \begin{equation}\label{TT-evo}
\xymatrix{
	   &                                                                      & \mathfrak{H} \mathcal{T}^{*}\mathcal{Q}    \times \mathbb{R} \ar@(ul,ur)^{ \mathcal{S}_{F}} \ar[ldd]_{\upsilon} \ar[rdd]^{\hat{\tau}_{\mathcal{T}^*\mathcal{Q}}}&&                                                                     \\  \\
											  & T\mathcal{Q} \ar[rdd]_{\tau^0_{\mathcal{Q}}}   && \mathcal{T}^*\mathcal{Q} \ar[ldd]^{\pi^0_{\mathcal{Q}}}
										\\	  \\&&\mathcal{Q}
										}
	\end{equation}
 where $\upsilon$ is the projection defined by $\upsilon= pr_{T\mathcal{Q}} \circ \mathcal{T}\pi_{\mathcal{Q}}^0 \circ j$. The local expression of this mapping is
 \begin{equation}\label{tilde-map}
\upsilon:\mathfrak{H} \mathcal{T}^{*}\mathcal{Q} \times \mathbb{R}\longrightarrow T\mathcal{Q},\qquad (q^i,p_i,z,\dot{q}^i,\dot{p}_{i}, u)\mapsto (q^i,\dot{q}^i).
 \end{equation}
 Here, $\hat{\tau}_{\mathcal{T}^*\mathcal{Q}}$ is the projection taking a pair $(V,u)$ mapping to $\tau_{\mathcal{T}^*\mathcal{Q}}(V)$, that is,
 \begin{equation}
 \hat{\tau}_{\mathcal{T}^*\mathcal{Q}}:\mathfrak{H} \mathcal{T}^{*}\mathcal{Q}    \times \mathbb{R}\longrightarrow \mathcal{T}^*\mathcal{Q},\qquad (q^i,p_i,z,\dot{q}^i,\dot{p}_{i}, u)\mapsto (q^i,p_i,z)
 \end{equation}
using that $ \mathfrak{H} \mathcal{T}^{*}\mathcal{Q}  $ is a vector subbundle of $ T\mathcal{T}^*\mathcal{Q} $.
 Consider now a (possibly non-horizontal) Lagrangian submanifold $\mathcal{S}_{F}$ of the symplectic manifold $\big(\mathfrak{H} \mathcal{T}^{*}\mathcal{Q}   \times \mathbb{R},\omega_\eta\big)$.
This submanifold projects to a submanifold $\upsilon(\mathcal{S}_{F})$ of $T\mathcal{Q}$ by the mapping $\upsilon$. Note that, $\upsilon(\mathcal{S}_{F})$  determines an implicit differential equation on the base $\mathcal{Q}$. The purpose of the Hamilton-Jacobi theory is to  retrieve solutions of $\mathcal{S}_{F}$, provided the solutions of the projected submanifold $\upsilon(\mathcal{S}_{F})$. 

In a similar way
as in (evolution) Hamilton-Jacobi Theorem \ref{HJT-Evo-Con-1}, in order to lift the solutions in $Q$ to $T^{*}Q$, we are still in need of a prolongation $\mathcal{T}^*W$ on $Q$, but two ingredients
of the theory are missing. One is that the base manifold, denoted by $C= \hat{\tau}_{\mathcal{T}^*\mathcal{Q}}(\mathcal{S}_{F})$, is not necessarily the whole $\mathcal{T}^*Q$, but possibly a proper submanifold of it.
The second is the nonexistence of a Hamiltonian vector field due to the
implicit character of the equations. 

We consider the restriction of $\mathcal{S}_{F}$ to the image  of $\mathcal{T}^*W$. This reads the following restricted submanifold 
\begin{equation} \label{S-Im-Evo}
\begin{split}
\mathcal{S}
_{F}\vert_{im(\mathcal{T}^*W)}&=\Big \{\left(q^i,\frac{\partial W}{\partial q^i}, W(q); \left.\frac{\partial F}{\partial p_i}\right\vert_{im(\mathcal{T}^*W)},\left.
-\frac{\partial F}{\partial q^i}\right\vert_{im(\mathcal{T}^*W)}-\left.\frac{\partial W}{\partial q^i}\frac{\partial F}{\partial z}\right\vert_{im(\mathcal{T}^*W)},\left.\frac{\partial F}{\partial z}\right\vert_{im(\mathcal{T}^*W)}\right)
\\ & \hspace{3.5cm} \in \mathfrak{H} \mathcal{T}^{*}\mathcal{Q}    \times \mathbb{R}:\left.\frac{\partial F}{\partial \lambda^a}\right\vert_{im(\mathcal{T}^*W)}=0\Big\}.
\end{split}
\end{equation}
If the Lagrangian submanifold $\mathcal{S}_{F}$ is the image of the pair $(\varepsilon_H,\mathcal{R}(H))$, then $\mathcal{S}_F\vert_{im(\mathcal{T}^*W)}$ reduces to the image   of $(\varepsilon_H,\mathcal{R}(H)) \circ \mathcal{T}^*W$. The submanifold $\mathcal{S}_F\vert_{im(\mathcal{T}^*W)}$ depends only on the base coordinates $(q^i)$. This enables us  to project it into  a submanifold $\mathcal{S}_F^{\mathcal{T}^*W}$ of $T\mathcal{Q}$ via the mapping $\upsilon$ as 
\begin{equation} \label{S-gamma}
\mathcal{S}_F^{\mathcal{T}^*W}=\upsilon \left( \mathcal{S}_F\vert_{im(\mathcal{T}^*W)}\right)=\left \{\left(q^i,\left.\frac{\partial F}{\partial p_i}\right\vert_{im(\mathcal{T}^*W)}\right )\in T\mathcal{Q}:\left.\frac{\partial F}{\partial \lambda^a}\right\vert_{im(\mathcal{T}^*W)}=0\right\}.
\end{equation}
In general, the submanifold $\mathcal{S}_F^{\mathcal{T}^*W}$ defines an implicit differential equation on $\mathcal{Q}$ given by
\begin{equation}\label{evo-eq-imp-proj}
\frac{dq^i}{dt}	=\left.\frac{\partial F}{\partial p_i}\right\vert_{im(\mathcal{T}^*W)},   \qquad \left.\frac{\partial F}{\partial \lambda^a}\right\vert_{im(\mathcal{T}^*W)}=0.
\end{equation}
We state the following Hamilton-Jacobi theorem for the implicit evolution Hamiltonian systems. 
\begin{theorem} \label{nHJT1-evo}
The following two conditions are equivalent for a real valued function $W=W(q)$ on $\mathcal{Q}$: 
\begin{enumerate}
\item The Lagrangian submanifold $\mathcal{S}_F$ in (\ref{S-Loc-Evo}) and the submanifold $\mathcal{S}_F^{\mathcal{T}^*W}$ in (\ref{S-gamma}) are $\mathcal{T}^*W$-related, that is
\begin{equation}
T\mathcal{T}^*W(\mathcal{S}_F^{\mathcal{T}^*W})=\mathcal{S}_F\vert_{im(\mathcal{T}^*W)}
\end{equation}
\item $dF(q^i,\partial W / \partial q^i, W,\lambda^a) =0$,
where $F$ is the Morse family generating $\mathcal{S}_F$.
\end{enumerate}
\end{theorem}
In terms of the local coordinates, the second condition in Theorem \ref{nHJT1-evo} establishes the implicit evolution Hamilton-Jacobi equation as
 \begin{equation} \label{Imp-HJ-loc}
F\left(q^i,\frac{\partial W}{\partial q^i},W(q),\lambda^\alpha\right)=c.
\end{equation}
where $c$ is the constant of integration.

\section{Example: The HJ Theory for Lagrangian Dynamics on Contact Manifolds}\label{Sec-Lag}

\subsection{Contact Lagrangian Dynamics}

Consider a Lagrangian function $L=L(q^i,\dot{q}^i,z)$ defined on the extended tangent bundle $\mathcal{T}\mathcal{Q}$. 
The Herglotz equations (also known as the generalized
Euler-Lagrange equations) are presented as the following system of equations \cite{de2020review,Guenther,Herglotz}:
\begin{equation}\label{Herglotz}
\frac{dq^i}{dt}=\dot{q}^i,\qquad 
\frac{\partial L}{\partial q^i} - \frac{d}{dt}\Big(\frac{\partial L}{\partial {\dot q}^i} \Big)
+ \frac{\partial L}{\partial z}\frac{\partial L}{\partial {\dot q}^i} = 0, \qquad  \frac{dz}{dt} = L(q^i,\dot{q}^i,z).
\end{equation}
By simply taking the fiber derivative of $L$ with respect to the velocity, one arrives the momenta, i.e., $p_i=\partial L / \partial \dot{q}^i$. A Lagrangian function is called regular (non-degenerate) if the rank of the Hessian matrix $[\partial^2 L /\partial \dot{q}^i \partial \dot{q}^j]$ is maximal. Recall that for regular Lagrangians, one can, at least locally, solve $\dot{q}^i$ in terms of $(q^i,p_i,z)$ and the Herglotz equations \eqref{Herglotz} can be easily (Legendre) transformed into the equations of contact Hamiltonian dynamics present in \eqref{conham}. See that for regular Lagrangians, \eqref{Herglotz} is a system of explicit differential equations. 

If a Lagrangian function is not regular, that is, referring to the equation $p_i=\partial L / \partial \dot{q}^i$ if one cannot (not even locally) solve $\dot{q}^i$ in terms of $(q^i,p_i,z)$, then the system \eqref{Herglotz} is a system of  implicit differential equations. For this case, let us rewrite \eqref{Herglotz} as 
\begin{equation}
\frac{dq^i}{dt}=\dot{q}^i, 
\qquad 
\frac{d{p}_i}{dt}
 = 
\frac{\partial L}{\partial q^i} 
+ \frac{\partial L}{\partial z}\frac{\partial L}{\partial {\dot q}^i} ,
\qquad 
\frac{dz}{dt} = L(q,\dot{q},z), \qquad p_i=\frac{\partial L} { \partial \dot{q}^i} .
\end{equation}
Consider now the Morse family
\begin{equation}\label{Morse-exp}
E(q,p,z,\dot{q})=\dot{q}^ip_i-L(q,\dot{q},z).
\end{equation} 
defined on the Whitney sum $\mathcal{T}^*\mathcal{Q}\times_{\mathcal{Q}\times \mathbb{R}}\mathcal{T}\mathcal{Q}$ over the base manifold $\mathcal{T}^*\mathcal{Q}$. Here, the velocity $\dot{q}^i$ is the Lagrange multiplier in the Morse family $E$. By substituting $E$ into the generic formulation in \eqref{N_E} one arrives at the Legendrian submanifold considering the Herglotz equations \eqref{Herglotz}, see \cite{esen2021contact} 
\begin{equation}\label{N_L} 
\mathcal{N}_L=\left\{\left(q^i,p_i,z, \dot{q}^i, \frac{\partial L}{\partial z}\frac{\partial L}{\partial \dot{q}^i}+\frac{\partial L}{\partial q^i},L, -\frac{\partial L}{\partial z}\right)\in \mathcal{T}\mathcal{T}^*\mathcal{Q} : ~ p_i=\frac{\partial L} { \partial \dot{q}^i} \right\}\subset \mathcal{T}\mathcal{T}^*\mathcal{Q}. 
\end{equation}

Let $\gamma(q,z)=(dW(q),z)$ be a section of the fibration $\mathcal{T}^*\mathcal{Q}\mapsto \mathcal{Q}   \times \mathbb{R}$ in the form of \eqref{gamma}. Here, $W$ is a real valued function on $\mathcal{Q}$. Consider the restriction of the Legendrian submanifold $\mathcal{N}_L$ to the image of $\gamma$. According to the formulation in \eqref{N-Im-E}, this gives 
\begin{equation} \label{N-Im-E-exp}  
\mathcal{N}_L\vert_{im(\gamma)} 
 =
\Big \{
\Big(q^i,\frac{\partial W}{\partial q^i},z; \dot{q}^i,\frac{\partial L}{\partial z}\frac{\partial L}{\partial \dot{q}^i}+\frac{\partial L}{\partial q^i},  
L,  - \frac{\partial L}{\partial z} \Big )
\in \mathcal{T}\mathcal{T}^*\mathcal{Q}  :\frac{\partial W}{\partial q^i}=\frac{\partial L} { \partial \dot{q}^i} \Big\}. 
\end{equation}
The submanifold $\mathcal{N}_L\vert_{im(\gamma)}$ depends only on the coordinates coordinates $(q^i,z,\dot{q}^i)$. So that, referring to \eqref{N-gamma}, we can project it to  a submanifold $\mathcal{N}_L^{\gamma}$ of $T(\mathcal{Q}\times \mathbb{R})$ as
\begin{equation} \label{N-gamma-} 
\mathcal{N}_L^{\gamma}= 
\Big \{
\Big(q^i,z; \dot{q}^i,  
L (q,\dot{q},z) \Big )
\in T(\mathcal{Q}\times \mathbb{R})  :\frac{\partial W}{\partial q^i}=\frac{\partial L} { \partial \dot{q}^i} \Big\} .
\end{equation}
The submanifold $\mathcal{N}_L^{\gamma}$ determines a system of implicit differential equations on $\mathcal{Q}\times \mathbb{R}$ given as
\begin{equation}\label{pro-con-equu}
\frac{dq^i}{dt}=\dot{q}^i, \qquad \frac{dz}{dt}=L (q,\dot{q},z) ,  \qquad \frac{\partial W}{\partial q^i}(q)=\frac{\partial L} { \partial \dot{q}^i}(q,\dot{q},z). 
 \end{equation}
Referring to Hamilton-Jacobi Theorem \ref{HJ-imp-con}, 
we are ready to state the Hamilton-Jacobi equation \eqref{HJ-imp-exp-cont} for non-regular Lagrangian dynamics as 
 \begin{equation}
     \frac{\partial W}{\partial q^i \partial q^j} \dot{q}^j - \frac{\partial L}{\partial q^i} - 
  \frac{\partial W}{\partial q^i}  \frac{\partial L}{\partial z}
  =0, \qquad 
 \frac{\partial W}{\partial q^i}= \frac{\partial L}{\partial \dot{q}^i}.
   \end{equation}
 A solution $W$ to this equations takes a solution of the projected system \eqref{pro-con-equu} to the Herglotz equations in \eqref{Herglotz}.

\subsection{Evolution Lagrangian Dynamics}

One can obtain the Lagrangian formalism for the evolution vector field by using a nonlinear nonholonomic action principle, as shown in \cite{simoes2020geometry}. The resulting equations are the evolution Herglotz equations
\begin{equation}\label{Herglotz-evo-} 
	\frac{dq^i}{dt}= \dot{q}^i,\qquad \frac{\partial L}{\partial q^i} - \frac{d}{dt}\Big(\frac{\partial L}{\partial {\dot q}^i} \Big)
	+ \frac{\partial L}{\partial z}\frac{\partial L}{\partial {\dot q}^i} = 0, \qquad 
	\frac{dz}{dt} = \dot{q}^i \frac{\partial L}{\partial \dot{q}^i}. 
\end{equation}
For regular Lagrangians, the evolution Herglotz equations \eqref{Herglotz-evo-} are transformed into the evolution Hamilton equations \eqref{evo-eq}. In this case, the dynamical equations are explicit. 
For non-regular Lagrangian functions, the system \eqref{Herglotz-evo-} is implicit. Accordingly, recalling the definition $p_i=\partial L / \partial \dot{q}^i$, we can rewrite \eqref{Herglotz-evo-} as 
\begin{equation}
\frac{dq^i}{dt}=\dot{q}^i, 
\qquad 
\frac{d{p}_i}{dt}
 = 
\frac{\partial L}{\partial q^i} 
+ \frac{\partial L}{\partial z}\frac{\partial L}{\partial {\dot q}^i} ,
\qquad 
\frac{dz}{dt} =  \dot{q}^i \frac{\partial L}{\partial \dot{q}^i}, \qquad p_i=\frac{\partial L} { \partial \dot{q}^i} .
\end{equation}
Referring to the Morse family $E=E(q,p,z,\dot{q})$ in \eqref{Morse-exp} 
defined on the Whitney sum $\mathcal{T}^*\mathcal{Q}\times_{\mathcal{Q}\times \mathbb{R}}\mathcal{T}\mathcal{Q}$ and to the generic formulation in \eqref{S-Loc-Evo-F} one arrives at the Lagrangian submanifold that determines the evolution Herglotz equations \eqref{Herglotz-evo-}, see \cite{esen2021contact}.
\begin{equation}\label{S_L-exp} 
\mathcal{S}_L=\left\{\left(q^i,p_i,z, \dot{q}^i, \frac{\partial L}{\partial z}\frac{\partial L}{\partial \dot{q}^i}+\frac{\partial L}{\partial q^i}, -\frac{\partial L}{\partial z}\right)\in \mathfrak{H} \mathcal{T}^{*}\mathcal{Q} \times \mathbb{R} : ~ p_i=\frac{\partial L} { \partial \dot{q}^i} \right\}.
\end{equation}

Let $W$ be a real valued function on $\mathcal{Q}$, we denote the restriction of $\mathcal{S}_L$ to the image space of the first prolongation $\mathcal{T}^*W$ as $\mathcal{S}
_{L}\vert_{im(\mathcal{T}^*W)}$. Referring to the formulation given in \eqref{S-Im-Evo}, we compute 
\begin{equation} \label{S-Im-Evo-exp}
\begin{split}
\mathcal{S}
_{L}\vert_{im(\mathcal{T}^*W)}&=\Big \{\big(q^i,\frac{\partial W}{\partial q^i}, W(q); \dot{q}^i,\left.\frac{\partial L}{\partial z}\frac{\partial L}{\partial \dot{q}^i}\right\vert_{im(\mathcal{T}^*W)}+\left.\frac{\partial L}{\partial q^i} \right\vert_{im(\mathcal{T}^*W)},\left.-\frac{\partial L}{\partial z}\right\vert_{im(\mathcal{T}^*W)}\big )
\\ & \hspace{3.5cm} \in \mathfrak{H} \mathcal{T}^{*}\mathcal{Q}    \times \mathbb{R}:\frac{\partial W}{\partial q^i}=\left.\frac{\partial L} { \partial \dot{q}^i}\right\vert_{im(\mathcal{T}^*W)}  \Big\}.
\end{split}
\end{equation}
According to \eqref{S-gamma}, we project $\mathcal{S}
_{L}\vert_{im(\mathcal{T}^*W)}$ to the tangent bundle as $T\mathcal{Q}$ as
\begin{equation} \label{S-gamma-exp}
\mathcal{S}_L^{\mathcal{T}^*W}= \left \{\left(q^i,\dot{q}^i \right ) \in T\mathcal{Q}:\frac{\partial W}{\partial q^i}=\left.\frac{\partial L} { \partial \dot{q}^i}\right\vert_{im(\mathcal{T}^*W)} \right\}.
\end{equation}
This corresponds to the implicit system of equations
\begin{equation}\label{reduced-evo-exp}
\frac{dq^i}{dt}=\dot{q}^i, 
\qquad 
\frac{\partial W}{\partial q^i}=\left.\frac{\partial L} { \partial \dot{q}^i}\right\vert_{im(\mathcal{T}^*W)}. 
\end{equation}
In general, the submanifold $\mathcal{S}_F^{\mathcal{T}^*W}$ defines an implicit differential equation on $Q$. The Hamilton-Jacobi Theorem \eqref{nHJT1-evo} gives the implicit Hamilton-Jacobi equation \eqref{Imp-HJ-loc} for this case as 
\begin{equation}
\dot{q}^i\frac{\partial W}{\partial q^i}-L(q,\dot{q},W(q))=cst.
 \end{equation}
A solution $W$ of this differential equation lifts a  solution of the projected dynamics \eqref{reduced-evo-exp} to a solution of the evolution Herglotz equations
\eqref{Herglotz-evo-}.

    \section{Conclusions}
    In this paper, on a contact manifold, we have introduced implicit contact Hamiltonian dynamics as a Legendrian submanifold of the tangent contact manifold. Two Hamilton-Jacobi theorems have been proposed: the first Theorem \ref{HJT-Con} provides a Hamilton-Jacobi theory only for explicit contact Hamiltonian dynamics,  whereas Theorem \ref{HJ-imp-con} proposes a Hamilton--Jacobi equation for implicit contact Hamiltonian dynamics as well. Evidently, Theorem \ref{HJ-imp-con} generalizes Theorem \ref{HJT-Con}.  Additionally, we have introduced implicit evolution Hamiltonian dynamics as a Lagrangian submanifold of a symplectic embedded space of the tangent contact manifold. A Hamilton-Jacobi equation has been stated for the evolutionary dynamics too. Theorem \ref{HJT-Evo-Con-1} provides a Hamilton-Jacobi theory for explicit evolution Hamiltonian dynamics, whilst Theorem \ref{nHJT1-evo} generalizes the Hamilton-Jacobi theory to implicit evolution dynamics. 
    
    \section{Acknowledgments}

We acknowledge the financial support from the Spanish Ministerio de Ciencia, Innovacion y Universidades Grant PID2019-106715GB-C21 and  the Severo Ochoa Programme for Centres of Excellence in R\&D" (CEX2019-000904-S).
Manuel Lainz wishes to thank MICINN and ICMAT for a FPI-Severo Ochoa predoctoral
contract PRE2018-083203.

\appendix

\newpage
\section{Morse Families and Special Symplectic Manifolds.} \label{appendix-1}

\textbf{Morse Families.} 
Let $\left(R,\tau, N \right) $ be a fiber bundle.
The vertical bundle $VR$ over $R$ is the space of vertical vectors $U\in TR$
satisfying $T\tau\left( U\right) =0$. The conormal bundle
of $VR$ is defined by
\begin{equation*}
V^{0}R=\left\{ \alpha\in T^{\ast}R:\left\langle \alpha,U\right\rangle
=0,\forall U\in VR\right\} .
\end{equation*}
Let $E$ be a real-valued function on $R$, then the image of its exterior derivative is a submanifold of $T^{\ast}R $. We say that $E$
is a Morse family (or an energy function) if
\begin{equation}
T_{z}im\left( dE\right) +T_{z}V^{0}R=T_zT^{\ast}R,  \label{MorseReq}
\end{equation}
for all $z\in im\left( dE\right) \cap V^{0}R$. A Morse family defined on $\left( R,\tau ,N\right) $ generates a Lagrangian submanifold of the canonical symplectic manifold $\left(
T^{\ast }N,\Omega_N\right)$ in the following way \cite{benenti2011hamiltonian,LiMa87,weinstein1981symplectic}. 
:
\begin{equation}
S_N=\left\{ w \in T^{\ast }N:T^{\ast }\tau (w
)=dE\left( z\right), \text{ for some } z\in T^*R \right\} . \label{LagSub}
\end{equation}%
 In this case, we say that $S_N$ is generated by the Morse family $E$. Note that, in the definition of $S_N$, there is an intrinsic requirement that $\tau \left( z\right)
=\pi _{N}\left( w \right) $. The inverse of this statement is also true, that is, any Lagrangian submanifold is generated by a Morse family. This is known as the Maslov-Hörmander theorem (or generalized Poincar\'{e} lemma). Assume that $N$ is equipped with local coordinates $(x^a)$, and consider the bundle local coordinates $(x^a,\lambda^\alpha)$ on the total space $R$. In this picture, a function $E$ is called a Morse family if the rank of the matrix%
\begin{equation}
\left( \frac{\partial ^{2}E}{\partial x^a x^b} \quad \frac{%
\partial ^{2}E}{\partial x^a \partial \lambda^\alpha }\right)  \label{MorseCon}
\end{equation}%
is maximal. In such a case, the Lagrangian submanifold (\ref{LagSub}) generated by $E$ locally looks like
\begin{equation} \label{MFGen}
S_N=\left \{\left(x^a,\frac{\partial E}{\partial x^a}(x,\lambda)\right )\in T^*N: \frac{\partial E}{\partial \lambda^\alpha}(x,\lambda)=0 \right \}.
\end{equation}
See that the dimension of $S_N$ is half of the dimension of $T^*N$, and that the canonical symplectic two-form $\Omega$ vanishes on $S_N$.

\textbf{Special Symplectic Structures.}  
 Let $P$ be a symplectic manifold carrying an exact symplectic two-form $%
\Omega=d\Theta$. Assume also that, $P$ is the total space of a fibre bundle $(P,\pi,N)$. A special symplectic structure is a quintuple $
(P,\pi,N,\Theta,\chi)$ where $%
\chi$ is a fiber preserving symplectic
diffeomorphism from $P$ to the cotangent bundle $T^{\ast}N$ \cite{LaSnTu75,SnTu73}. Here, $\chi$ can uniquely be characterized by
\begin{equation} \label{chi-}
\left\langle \chi(p),\pi_{\ast}X(n)\right\rangle =\left\langle
\Theta(p),X(p)\right\rangle
\end{equation}
for a vector field $X$ on $P$, for any point $p$ in $P$ where $\pi(p)=n$. Note that, the pairing on the left hand side of (\ref{chi-}) is the natural pairing between the cotangent space $T^*_nN$ and the tangent space $T_nN$. The pairing on the right hand side of (\ref{chi-}) is the one between the cotangent space $T^*_pP$ and the tangent space $T_pP$. Here is a diagram exhibiting the special symplectic structure.
\begin{equation} \label{sss}
\xymatrix{T^{\ast }N \ar[d]_{\pi_{N}} &&P
\ar[dll]^{\pi} \ar[ll]_{\chi}
\\N }
\end{equation}
The pair $(P,\Omega)$ is called the underlying symplectic manifold of the special symplectic structure $%
(P,\pi,N,\Theta,\chi)$. 
\bigskip

\noindent 
Let $
(P,\pi,N,\Theta,\chi)$ be a special symplectic structure. Assume also that $S_P$ is a Lagrangian submanifold of $P$. The image $S_N=\chi(S_P)$ of $S_P$ is a Lagrangian submanifold of $T^{\ast}N$. By referring to the generalized Poincar\'{e} lemma presented in the previous subsection, we argue that the Lagrangian submanifold $S_N$ can locally be generated by a Morse family $E$ on a fiber bundle $(R,\tau,N)$. Accordingly, we call the Morse family $E$ a generator of both $S_N$ and $S_P$ since they are the same up to $\chi$. The following diagram summarizes this discussion by equipping a Morse family  to a special symplectic structure (\ref{sss}).
 \begin{equation} \label{Morse-Gen}
 \xymatrix{
\mathbb{R}& R \ar[d]^{\tau}\ar[l]_{E}& T^*N \ar@(ul,ur)^{S_N}  \ar[d]_{\pi_N}& &P \ar@(ul,ur)^{S_P}   \ar[ll]_{\chi} \ar[dll]^{\pi}\\ &
N \ar@{=}[r]& N
}
\end{equation}

\bibliographystyle{amsplain}
\bibliography{references}
\end{document}